\documentclass[a4paper,fleqn]{cas-sc}

\usepackage[square,sort,comma,compress,numbers]{natbib}
\usepackage{graphicx}
\usepackage{graphics}
\usepackage[section]{placeins}
\usepackage{caption}
\usepackage{subcaption}
\usepackage{epsfig}
\usepackage{booktabs}
\usepackage{color}
\usepackage{amsmath}
\usepackage{mathtools}
\graphicspath{ {./figures/} }
\usepackage{booktabs}
\usepackage{multirow}
\usepackage{pifont}
\usepackage{inputenc}
\usepackage{tikz}
\usetikzlibrary{shapes,arrows}
\usetikzlibrary{shapes.geometric}
\usetikzlibrary{calc}

\newcommand{\ary}[1]{\boldsymbol{\mathsf{#1}}}
\newcolumntype{L}{>{$}l<{$}}
\newcolumntype{C}{>{$}c<{$}}
\usepackage{lineno}
\usepackage{mathrsfs}
\usepackage{color}
\usepackage{xcolor}
\def\tsc#1{\csdef{#1}{\textsc{\lowercase{#1}}\xspace}}
\tsc{WGM}
\tsc{QE}
\tsc{EP}
\tsc{PMS}
\tsc{BEC}
\tsc{DE}

\begin{document}
\let\WriteBookmarks\relax
\def\floatpagepagefraction{1}
\def\textpagefraction{.001}
\shorttitle{ }
\shortauthors{\today}

\title [mode = title]{Computational bifurcation analysis of hyperelastic thin shells}

\author[1]{Zhaowei Liu}[orcid=0000-0002-0572-7415]

\address[1]{College of Mechanics and Materials, Hohai University, Nanjing, 211100, China}
\address[2]{Glasgow Computational Engineering Centre, James Watt School of Engineering, University of Glasgow, Glasgow, G12 8LT, United Kingdom}

\author[2]{Andrew McBride}[orcid = 0000-0001-7153-3777]

\author[2]{Abhishek Ghosh}[orcid = 0000-0003-3066-4492]

\author[3]{Luca Heltai}[orcid = 0000-0001-5514-4683]
\address[3]{SISSA (International School for Advanced Studies), Via Bonomea 265, 34136 Trieste, Italy}

\author[4] {Weicheng Huang}[]
\address[4] {School of Mechanical Engineering, Southeast University, Nanjing, 211189, China}

\author[1]{Tiantang Yu}[]

\author[2,5]{Paul Steinmann}[orcid = 0000-0003-1490-947X]

\address[5]{Institute of Applied Mechanics, Friedrich-Alexander Universit\"at Erlangen-N\"urnberg, D-91052, Erlangen, Germany}

\author[2]{Prashant Saxena}[orcid = 0000-0001-5071-726X]
\cormark[1] 
\ead{Prashant.Saxena@glasgow.ac.uk}

\cortext[cor1]{Corresponding author}

\begin{abstract}
The inflation of hyperelastic thin shells is an important and highly nonlinear problem that arises in multiple engineering applications involving severe kinematic and constitutive nonlinearities in addition to various instabilities.
We present an isogeometric approach to compute the inflation of hyperelastic thin shells, following the Kirchhoff-Love hypothesis and associated large deformation. 
Both the geometry and the deformation field are discretized using Catmull-Clark subdivision bases which provide the $C^1$-continuous finite element framework required for the Kirchhoff-Love shell formulation. 
To follow the complex nonlinear response of hyperelastic thin shells, the inflation is simulated incrementally, and each incremental step is solved via the Newton-Raphson method enriched with arc-length control. 
Eigenvalue analysis of the linear system after each incremental step allows for inducing bifurcation to a lower energy mode in case stability of the equilibrium is lost.
The proposed method is first validated using benchmarks, and then applied to engineering applications, where we demonstrate the ability to simulate large deformation and associated complex instabilities.
\end{abstract} 

\begin{highlights}
\item An accurate, efficient and robust isogeometric approach has been developed to simulate large inflation of nonlinear elastic thin shells.
\item The nonlinear thin shell formulation is implemented based on Catmull-Clark subdivision surfaces to analyse incompressible hyperelastic thin shells with arbitrary geometry.
\item {The proposed algorithm predicts the complex mechanical behaviour of highly nonlinear  inflated thin shells. It provides investigation tools for complex behaviours of such structures in soft robotics and actuators.}
\item Three distinct types of instabilities, namely snap-through buckling, loss of symmetry bifurcation and wrinkling are successfully investigated with the proposed method.
\end{highlights}

\maketitle

\section{Introduction}
\label{intro}

{ 
An accurate, efficient, and robust isogeometric approach is developed to simulate various instabilities that may occur during the large deformation of nonlinearly elastic thin shells. 
The method is capable of predicting the mechanical behaviour of  nonlinear thin shells with arbitrary geometry, and provides guidance for the design and application of highly deformable thin structures. 
}

Many natural rubbers, synthetic elastomers, and soft biological tissues can undergo large reversible deformations and demonstrate a nonlinear stress-strain relationship.
Hyperelastic constitutive laws are used to accurately model their nonlinear behaviour.
Moreover, most elastomers also demonstrate isotropy and incompressibility which can be used to simplify the constitutive relations.
Among the various choices of hyperelastic constitutive models available in the literature (neo-Hookean, Mooney--Rivlin, Ogden, Gent, Arruda--Boyce, to name a few),  the Mooney--Rivlin model is adopted here on account of its simplicity and its ability to capture the appropriate nonlinear response.
This model requires only two material constants along with the first and second invariants of the right Cauchy--Green deformation tensor to define the strain energy density function. 
It reduces to the neo-Hookean constitutive model when the material constant associated with the second invariant vanishes.


A particularly relevant and challenging application of these nonlinear materials emerges when one considers slender structures such as rods, membranes, plates, and shells.
In these structures, one or more characteristic dimensions are negligible with respect to the others.
They can be modelled as lower dimensional manifolds embedded in the three-dimensional space with appropriate kinematic simplifications.

In this work we consider the case of incompressible and hyperelastic thin shells, which can be regarded as two-dimensional surfaces with a small thickness in three dimensional space (membranes and plates are a particular subcase of these types of structures). Thin shells can undergo large deformation, even when subjected to small external loads, and their analysis  requires a careful consideration of kinematic, constitutive, and geometric nonlinearities.
The simplest kinematic simplification one can assume for shell models is based on the Kirchhoff hypothesis: lines perpendicular to the mid-surface remain straight and perpendicular to the mid-surface after deformation, neglecting out-of-plane shearing, and leading to the well established Kirchhoff-Love thin shell theory, which is governed by a fourth-order partial differential equation. The weak formulation of such equations seeks for solutions in spaces where the Hessians have a bounded $L^2$ norm, and requires a discretization where solutions are at least globally $C^1$ continuous~\cite{zienkiewicz2005finite}. Traditional finite element methods based on Lagrange polynomials only guarantee  $C^0$ continuity between elements, and one requires to resort to more exotic finite element approximations, such as Argyris finite element spaces on triangles~\cite{Argyris1968} or their equivalent on quadrilaterals~\cite{Kapl2021}. \citet{noels2008new} also proposed a discontinuous Galerkin  (DG) method for Kirchhoff-Love shell and meshless method~\cite{krysl1996analysis} is also able to overcome this problem.

An alternative approach that guarantees $C^1$ continuity comes from the isogeometric paradigm, where the basis functions of the finite element spaces are inspired by computer aided design (CAD) principles.
One of the first occurrences of this approach was presented by~\citet{cirakortiz2000} who developed a finite element approach based on Loop subdivision surfaces for Kirchhoff-love thin shells and later extended it to hyperelastic thin shells~\cite{Cirak:2001aa}.
Subdivision surfaces are a mature CAD tool that is widely used in the animation industry and engineering design. 
One of their strengths is that they guarantee global $C^1$ continuity for arbitrary control point topologies (including extraordinary vertices, i.e., vertices on a surface shared by a number of cells different from four, and hanging node vertices), and -- at the same time -- allow for fast evaluation using standard cubic spline functions in all ordinary patches (i.e., all patches where vertices are shared exactly by four cells). Shell formulations based on subdivision surfaces have been extended to applications including shell fracture~\cite{cirak2005cohesive}, shape optimization~\cite{BANDARA201862,chen2020acoustic}, fluid-structure interaction~\cite{cirak2007large}, structural-acoustic analysis~\cite{liu2018isogeometric,chen2022sample,chen2022multi} and piezoelectric shells~\cite{liu2022vibration}. 

The idea behind the use of more general non-uniform rational B-splines (NURBS) basis functions in the finite element context was proposed by  \citet{hughes2005isogeometric} in 2005 with application to elastic thin shells by~\citet{kiendl2009isogeometric}. The isogeometric approach was used by \citet{kiendl2015isogeometric} to analyse both compressible and incompressible hyperelastic thin shells, and extended to fluid structure interaction, coupling isogeometric BEM formulations~\cite{heltai2014nonsingular} with elastic thin shells in \cite{heltai2017natural}. \citet{takizawa2019isogeometric} derived a hyperelastic thin shell formulation with isogeometric discretization, where the out-of-plane deformation mapping is taken into account. \citet{tepole2015isogeometric} developed an isogeometric formulation of Kirchhoff-Love shells to analyse biological membranes. \citet{roohbakhshan2017efficient} also used an isogeometric rotation-free shell formulation to model soft tissues.  \citet{huynh2020elasto} studied elasto-plastic large deformation behavior of thin shell structures using the isogeometric approach. IGA is also extensively applied to analyse hyperelastic solids~\cite{bazilevs2008isogeometric,bernal2013isogeometric,hassani2015solution,du2020nitsche}.  The buckling of thin structures is also studied using IGA shell formulations. \citet{guo2019isogeometric} proposed an isogeometric analysis framework for the buckling analysis of trimmed elastic shells. \citet{verhelst2021stretch} presented a formulation of stretch-based material models for isogeometric Kirchhoff-Love shells and it was used to simulate tension wrinkling of a thin sheet. { An attractive feature of subdivision surfaces is that they can be evaluated using spline functions, while retaining a simple polygonal mesh data structure able to represent complex geometries and also allowing extraordinary vertices in the mesh, which enables local refinement and patch conforming, both challenges faced by NURBS.}

The inflation of thin structures made from rubber-like materials has numerous engineering applications, including tyres, airbags, air springs, buffers, pneumatic actuators~\cite{galley2019pneumatic}, and soft grippers~\cite{hao2017modeling}. The large deformations of inflated hyperelastic circular plates are well-studied with semi-analytical approaches~\cite{adkins1952large,hart1967large,yang1970axisymmetrical, SAXENA2019250}. 
The inflation of other axisymmetric thin structures, including cylindrical~\cite{khayat1992inflation,guo2001large,pamplona2006finite, REDDY2018203}, spherical~\cite{akkas1978dynamic,verron1999dynamic, XIE2016182} and toroidal~\cite{tamadapu2013finite,pamulaparthi2019instabilities, REDDY2017248} membranes, has also been investigated semi-analytically. 
\citet{holzapfel1996large} presented a general formulation of thin incompressible membranes to investigate biological tissues using the finite element method. \citet{bonet2000finite} analysed hyperelastic membranes containing an enclosed fluid with a finite element approach. \citet{rumpel2005efficient} also developed a finite element model for gas and fluid supported membrane and shell structures. 
The inflation of hyperelastic membranes and thin shells often involve in different types of instabilities. The most common instability phenomenon of inflating membranes is snap-through buckling. Hyperelastic membranes appear to lose stiffness at a critical point and suddenly undergo very large inflation with only a small pressure increment. This phenomenon is well studied~\cite{benedict1979determination,carroll1987pressure,khayat1992inflation,muller2002inflating,tamadapu2013finite}. Another well-known instability of hyperelastic membranes and thin shells is global buckling which breaks geometrical symmetries. It manifests as a bifurcation from symmetric deformation to asymmetric deformation in an inflated structure at a critical loading point. \citet{koiter1967stability} pioneered a method to evaluate the critical point of such instability by checking the sign of the second variation of the total potential energy. \citet{pamulaparthi2019instabilities} adopted this method to analyse the bifurcation  of inflated hyperelastic toroidal membranes.
For the nonlinear finite element approach, the stability can be checked by computing the eigenvalues of the stiffness matrix for each load step. The eigenvector corresponding to the zero eigenvalue indicates the direction of a possible bifurcation~\cite{de1988bifurcations}.

Most of the aforementioned instability analyses of hyperelastic thin structures are limited to simple geometries, and computations performed until the onset of bifurcation.
Rapid developments in manufacturing, soft-robotics, and biomedical engineering motivate the need to overcome these restrictions. 
To this end, an incompressible nonlinear hyperelastic thin shell formulation is implemented using an isogeometric approach based on Catmull-Clark subdivision surfaces.
The proposed method is able to handle both kinematic and constitutive nonlinearities of hyperelastic thin shells with arbitrary geometry.
We demonstrate that multiple types of instabilities associated with the large deformation of hyperelastic thin shells can be captured.
Furthermore, we are able to compute the important post-bifurcation response.
Specific attention is paid here to thin shells undergoing large inflation and the various resulting instabilities including limit point, wrinkling, and loss of symmetry bifurcations.




This manuscript is organized as follows. Section~\ref{sec:notations} introduces the notation and defines coordinate systems used throughout the manuscript. Section~\ref{sec:shell} introduced the thin shell formulation, where Section~\ref{sec:shell_geo} and~\ref{sec:shell_kin} define the geometric definitions and kinematics of nonlinear Kirchhoff-Love thin shells, respectively. Section~\ref{sec:constitutive} presents the constitutive relations of incompressible Mooney-Rivlin material adapted to the specific shell formulation. Thereafter, Section~\ref{sec:virtual_work} derives the governing equation of a nonlinear Kirchhoff-Love thin shell using the principle of virtual work. Section~\ref{sec:implementation} introduces the Catmull-Clark subdivision surfaces and discusses the implementation details of 
the isogeometric nonlinear finite element method for hyperelastic thin shell. In Section~\ref{sec:investigations}, the algorithm to simulate the nonlinear deformation of hyperelastic thin shells is illustrated, and the difficulties in capturing the snap-through and bifurcation bucklings of inflated hyperelastic thin structures are emphasised. Finally, four numerical examples are demonstrated: a circular plate problem validates the proposed method, and the snap-through buckling of a spherical shell is simulated with the results well agreeing with analytical solutions. The third numerical example investigates the bifurcation of a toroidal thin shell, while the last example simulates the folding and wrinkling of an inflated airbag.
\section{Notation}
\label{sec:notations}
\paragraph{Brackets:}
Square brackets $[ \, ]$ are used to group algebraic expressions. Round brackets $( \,)$ are used to denote the dependencies of a function.  If brackets are used to denote an interval, then $(\,)$ stands for an open interval and $[\,]$ is a closed interval. Curly brackets $\{\,\}$ are used to define sets.
\paragraph{Symbols:}
A variable typeset in a normal weight font represents a scalar. A bold weight font denotes a first or second-order tensor. An overline indicates that the variable is defined with respect to the reference configuration and if absent, the variable is defined with respect to the deformed configuration. A scalar variable with superscript or subscript indices normally represents the components of a vector or second-order tensor. Upright font is used to denote matrices and vectors.

Indices $i,j,k,\dots$ vary from $1$ to $3$ while $\alpha, \beta, \gamma,\dots$, used to indicate surface variable components, vary from $1$ to $2$. Einstein summation convention is used throughout.

The comma symbol in a subscript represents a partial derivative, for example, $A_{,\beta}$ is the partial derivative of $A$ with respect to the $\beta^{\text{th}}$ coordinate.
\paragraph{Coordinates:}
$ \mathbf{c}_i$ represent the basis vectors of an orthonormal system in three-dimensional Euclidean space and $x,y$ and $z$ are its coordinates. $\boldsymbol{\theta}^i$ denote the basis vectors in the local element space and $\theta^1,\theta^2$ and $\theta^3$ are its coordinates. The three covariant basis vectors for a surface point are denoted as $\mathbf a_i$, where $\mathbf a_1, \mathbf a_2$ are tangential vectors and $\mathbf a_3$ is the normal vector. 

\begin{nolinenumbers}
\section{Nonlinear Kirchhoff-Love shell formulation}
\label{sec:shell}
\subsection{Geometry}
\label{sec:shell_geo}
 Consider a shell in its reference configuration occupying a physical domain $\bar\Omega \subset \mathbb{R}^3$ as shown in \autoref{fig:shell_definition}. 
 Each point $\bar{\mathbf{r}} \in \bar\Omega$ is mapped from the parametric domain defined by the coordinate system $\{ \theta^1, \theta^2, \theta^3\}$.
 The Kirchhoff-Love hypothesis states that lines perpendicular to the mid-surface of the thin shell remain straight and perpendicular to the mid-surface after deformation.
 Hence, assuming the shell has a uniform thickness $\bar h$ in the reference configuration, the point $\bar{\mathbf r}$ in the shell-space can be defined using a point on the mid-surface $\bar\Gamma$ , denoted $\bar{\mathbf x} \in \bar\Gamma$, and the associated unit normal vector $\bar{\mathbf n}$ as
 \begin{equation}
\bar {\mathbf r} = \bar{\mathbf x} + \theta^3 \bar{\mathbf n},
\label{eq: r x n relation}
 \end{equation}
 where $\theta^3 \in [-{\bar h}/{2}, {\bar h}/{2}]$.
\begin{figure}[h]
\centering
  \includegraphics[width=0.5\linewidth]{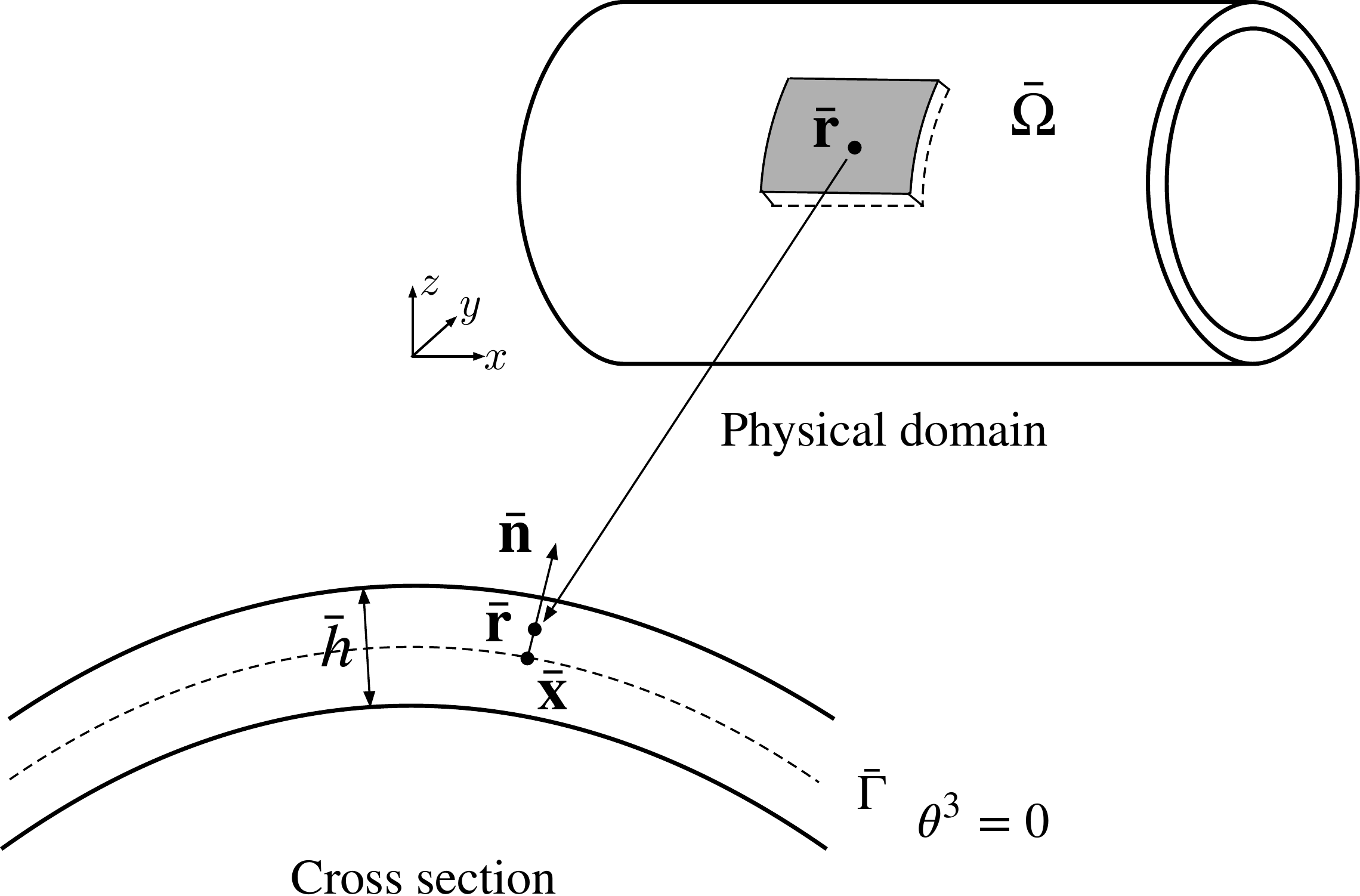}
\caption{A Kirchhoff--Love shell occupying a physical domain $\bar\Omega$. Each point $\bar{\mathbf{r}} \in \bar\Omega$ can be defined using quantities on the mid-surface $\bar\Gamma$ on the shell using $\bar {\mathbf r} = \bar{\mathbf x} + \theta^3 \bar{\mathbf n}$.}
\label{fig:shell_definition}
\end{figure}

The reference and deformed configurations must be carefully distinguished in order to study the finite deformations of hyperelastic shells. \autoref{fig:shell_deformed_reference} shows the reference and the deformed configuration of the mid-surface. Both configurations can be mapped from the parametric domain of the mid-surface. The points on the mid-surface in the reference and the deformed configurations are denoted by $\bar{\mathbf x}$ and $\mathbf x$, respectively. 
The mid-surface point in the deformed configuration ${\mathbf x}$ can be related to the mid-surface point in the reference configuration as
\begin{equation}
\mathbf x = \bar{\mathbf x} + \mathbf u,
\label{eq:reference_to_deformed}
\end{equation}
where $\mathbf u$ denotes the displacement. Moreover, the covariant basis vectors in the mid-surface of reference and the deformed configuration are computed as
\begin{equation}
\bar{\mathbf a}_\alpha = \frac{\partial \bar{\mathbf x}}{\partial \theta^\alpha}, 
\quad \text{and} \quad
{\mathbf a}_\alpha = \frac{\partial {\mathbf x}}{ \partial \theta^\alpha}.
\end{equation}
Thus, the unit normal vectors in the two configurations are defined by
\begin{equation}
\bar{\mathbf n} = \bar{\mathbf a}^3 = \frac{\bar{\mathbf a}_1 \times \bar{\mathbf a}_2}{\bar{J}} \equiv \bar{\mathbf a}_3, \quad \text{and} \quad  {\mathbf n} = {\mathbf a}^3 = \frac{{\mathbf a}_1 \times {\mathbf a}_2}{{J}} \equiv {\mathbf a}_3,
\end{equation}
and $\bar J$ and $J$ are the respective mid-surface Jacobians given by
\begin{equation}
\bar{J} = |\bar{\mathbf a}_1 \times \bar{\mathbf a}_2|,\quad \text{and} \quad J = | {\mathbf a}_1 \times {\mathbf a}_2|.
\end{equation}

\begin{figure}
\centering
  \includegraphics[width=0.6\linewidth]{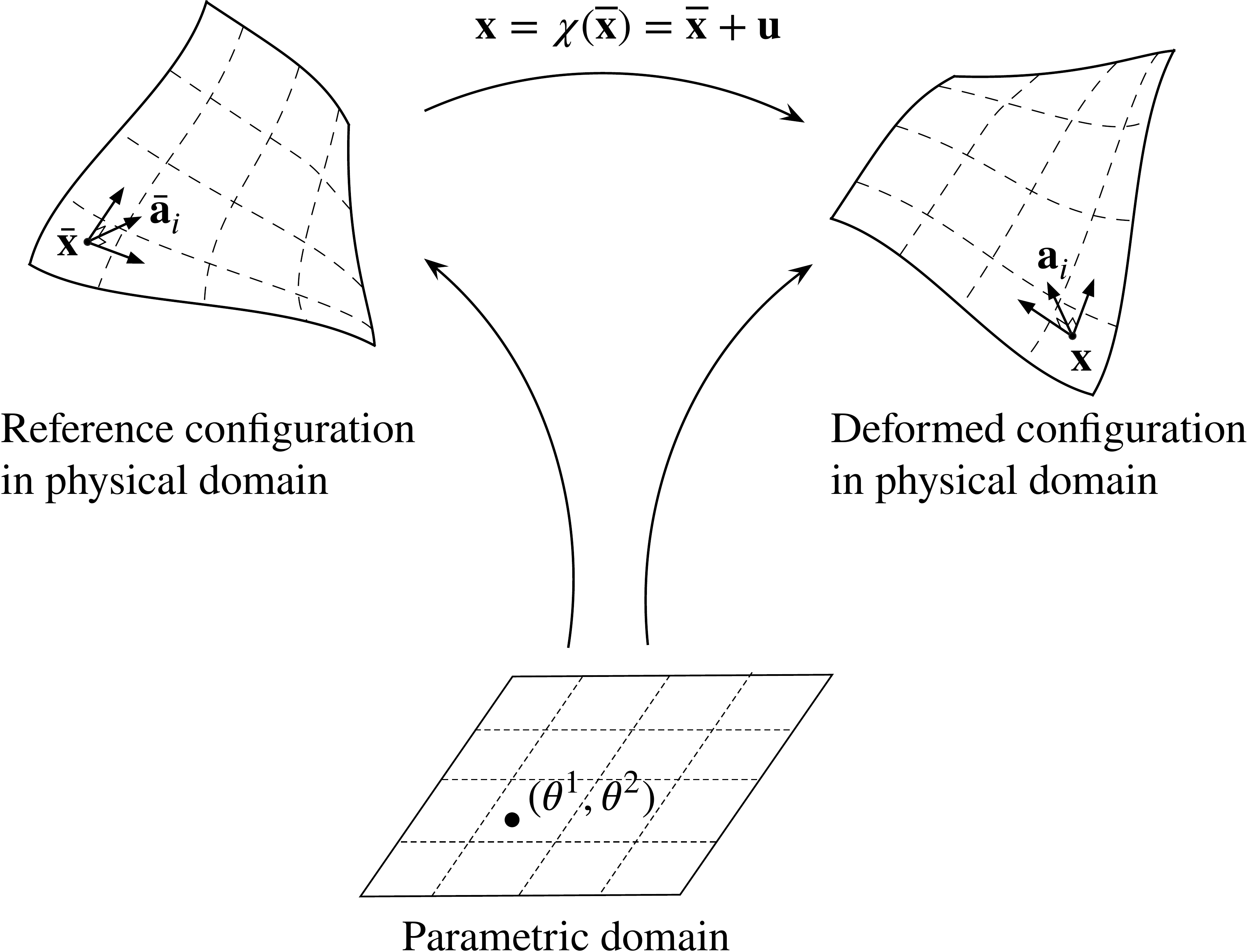}
\caption{Position vectors in the reference configuration $(\bar{\mathbf{x}} \in \bar\Gamma)$ of the mid-surface and deformed configuration $(\mathbf{x} \in \Gamma)$ of the mid-surface are related by the displacement vector $\mathbf{u}$.
$\bar\Gamma$ is spanned by the covariant basis $\{ \bar{\mathbf{a}}_1, \bar{\mathbf{a}}_2 \}$ while $\Gamma$ is spanned by the covariant basis $\{ \mathbf{a}_1, \mathbf{a}_2\}$.
}
\label{fig:shell_deformed_reference}
\end{figure}

The covariant components of the metric tensors for the mid-surface points $\bar{\mathbf x}$ and $\mathbf x$ are respectively given by
\begin{equation}
\bar{a}_{ij} = \bar{\mathbf a}_i \cdot \bar{\mathbf a}_j, \quad \text{and} \quad {a}_{ij} = {\mathbf a}_i \cdot {\mathbf a}_j.
\end{equation}
The corresponding contravariant metric tensors $\bar{a}^{ik}$ and $a^{ik}$ are defined by
\begin{equation}
\bar{a}^{ik}\bar{a}_{kj} = \delta^{i}_{j}, \quad \text{and} \quad a^{ik}a_{kj} = \delta^{i}_{j},
\label{eq:co_and_contra_metric}
\end{equation}
where $\delta^{i}_{j}$ denotes the Kronecker delta. 

The thickness stretch $\lambda_3$ for a finitely deformed shell is defined by
\begin{equation}
\lambda_3 = \frac{h}{\bar h},
\end{equation}
where $h = h(\theta^1,\theta^2)$ is the shell thickness after deformation.
We introduce a vector $\mathbf d$ to combine the thickness stretch and normal vector as
\begin{equation}
    \mathbf d = \lambda_3 \mathbf{a}_3,
\end{equation}
to write the position vector $\mathbf r$ of a point in the deformed configuration of the shell-space as
\begin{equation}
     {\mathbf r} = {\mathbf x} + \theta^3 {\mathbf d}, \quad \text{with} \quad \mathbf x = \mathbf x(\theta^1,\theta^2),\, \mathbf d = \mathbf d(\theta^1,\theta^2).
\end{equation}
Thus, the three-dimensional covariant basis vectors in the shell-space of the reference and the deformed configurations respectively follow as
\begin{equation}
\bar{\mathbf g}_\alpha = \frac{\partial \bar{\mathbf r}}{\partial \theta^\alpha} = \bar{\mathbf a}_\alpha + \theta^3 \bar{\mathbf a}_{3,\alpha},  \quad \bar{\mathbf g}_3 = \frac{\partial \bar{\mathbf r}}{\partial \theta^3} = \bar{\mathbf a}_3,
\label{eq:cov_tensors_r}
\end{equation}
and 
\begin{equation}
{\mathbf g}_\alpha = \frac{\partial {\mathbf r}}{\partial \theta^\alpha} = {\mathbf a}_\alpha + \theta^3 \mathbf{d}_{,\alpha}, \quad 
{\mathbf g}_3 = \frac{\partial {\mathbf r}}{\partial \theta^3} = {\mathbf d}.
\label{eq:cov_tensors_d}
\end{equation}
The components of the covariant metric tensors in the shell-space are given by
\begin{equation}
\bar{g}_{ij} = \bar{\mathbf g}_i \cdot \bar{\mathbf g}_j\quad \text{and} \quad g_{ij} = \mathbf g_i \cdot \mathbf g_j,
\end{equation}
and the contravariant components of the metric tensor at point $\mathbf r$ follow as
\begin{equation}
    \bar{g}^{ij} = \bar{\mathbf g}^{i} \cdot \bar{\mathbf g}^{j} \quad \text{and} \quad {g}^{ij} = {\mathbf g}^{i} \cdot {\mathbf g}^{j},
    \label{eq:contra_metric}
\end{equation}
where $\bar{\mathbf{g}}^i$ and ${\mathbf{g}}^i$ denotes the contravariant basis vectors in reference and deformed configuration of the shell-space defined by
\begin{equation}
\bar{\mathbf{g}}^i \cdot \bar{\mathbf{g}}_j = \delta^i_j \quad \text{and} \quad {\mathbf{g}}^i \cdot {\mathbf{g}}_j = \delta^i_j.
\label{eq:deform_tensor}
\end{equation}
\subsection{Kinematics}
\label{sec:shell_kin}
The three-dimensional deformation gradient in the shell-space is defined by 
\begin{equation}
\mathbf{F}=\mathbf{g}_{i} \otimes \bar{\mathbf{g}}^{i} .
\label{eq:deformation_gradient}
\end{equation}
The corresponding right Cauchy-Green deformation tensor follows as
\begin{equation}
\mathbf C = \mathbf F^{\mathrm T} \cdot \mathbf F = g_{ij}\bar{\mathbf{g}}^{i} \otimes \bar{\mathbf{g}}^{j},
\label{eq:Cauchy-Green_2}
\end{equation}
and its inverse is defined as
\begin{equation}
    {\mathbf{C}}^{-1} = \mathbf F^{-1} \cdot \mathbf F^{-\mathrm{T}} = {g}^{ij} \bar{\mathbf{g}}_i \otimes \bar{\mathbf{g}}_j .
    \label{eq:Cauchy-Green_inv}
\end{equation}
Thus, the components of ${\mathbf{C}}^{-1}$ are the components of the contravariant metric tensor defined in Equation~\eqref{eq:contra_metric}.
Upon substituting Equations~\eqref{eq:cov_tensors_d} into~\eqref{eq:Cauchy-Green_2}, the right Cauchy-Green deformation tensor can be expanded as
\begin{equation}
    \mathbf C=g_{\alpha \beta} \, \bar{\mathbf{g}}^{\alpha} \otimes \bar{\mathbf{g}}^{\beta}+\lambda_3^{2} \bar{\mathbf{g}}^{3} \otimes \bar {\mathbf{g}}^{3} + [\mathbf a_\alpha \cdot \mathbf d + \theta^3\mathbf d_{,\alpha}\cdot \mathbf d  ] \bar{\mathbf{g}}^{\alpha} \otimes \bar{\mathbf{g}}^{3} + [\mathbf a_\beta \cdot \mathbf d + \theta^3\mathbf d_{,\beta}\cdot \mathbf d ] \bar{\mathbf{g}}^{3} \otimes \bar{\mathbf{g}}^{\beta}.
\end{equation}
The third and fourth terms in the above equation are components corresponding to shear along the thickness direction. 
Since $\mathbf a_\alpha$ and $\mathbf n$ are always perpendicular to each other, $\mathbf a_\alpha \cdot \mathbf d = 0$. 
The term $\mathbf d_{,\alpha} \cdot \mathbf d = \lambda_3 \lambda_{3,\alpha} $
can be neglected on account of the long-wave assumption~\cite{maurin2016wave} for the change of the thickness stretch with mid-surface coordinates $\lambda_{3,\alpha} \approx 0$.
Therefore, there is no out-of-plane shear strain component in the right Cauchy-Green deformation tensor.
This leads to the following simplified form
\begin{equation}
\mathbf C=g_{\alpha \beta} \, \bar{\mathbf{g}}^{\alpha} \otimes \bar{\mathbf{g}}^{\beta}+\lambda_3^{2} \bar{\mathbf{g}}^{3} \otimes \bar {\mathbf{g}}^{3}.
\label{eq:Cauchy-Green_3}
\end{equation}
Ignoring the high order terms, the components of the covariant metric tensors read
\begin{equation}
g_{\alpha \beta}=a_{\alpha \beta}-2 \theta^{3} b_{\alpha \beta},
\label{eq:covariant_metric_tensor}
\end{equation}
with the first and second fundamental forms of the mid-surface of the deformed configuration
\begin{equation}
a_{\alpha\beta} = \mathbf{a}_\alpha \cdot \mathbf{a}_\beta\quad \text{and} \quad b_{\alpha\beta} = \mathbf{a}_{\alpha,\beta} \cdot \mathbf d =\lambda_3 \mathbf{a}_{\alpha,\beta} \cdot \mathbf a_3.
\end{equation}
Finally, the Green-Lagrange strain tensor can be expressed as
\begin{equation}
\mathbf E =\frac{1}{2}[\mathbf C - \mathbf I] = \frac{1}{2} [g_{\alpha \beta} - \bar g_{\alpha \beta}] \bar{\mathbf{g}}^{\alpha} \otimes \bar{\mathbf{g}}^{\beta}+ \frac{1}{2} [\lambda_3^{2}-1] \bar{\mathbf{g}}^{3} \otimes \bar{\mathbf{g}}^{3}.
\end{equation}
Similar to Equation~\eqref{eq:covariant_metric_tensor}, the components of the covariant metric tensor in the reference configuration read
\begin{equation}
\bar{g}_{\alpha \beta}=\bar{a}_{\alpha \beta}-2 \theta^{3} \bar{b}_{\alpha \beta},
\end{equation}
with the first and second fundamental forms of the mid-surface of the reference configuration
\begin{equation}
\bar{a}_{\alpha\beta} = \bar{\mathbf{a}}_\alpha \cdot \bar{\mathbf{a}}_\beta\quad \text{and} \quad \bar{b}_{\alpha\beta} = \bar{\mathbf{a}}_{\alpha,\beta} \cdot \bar{\mathbf a}_3.
\end{equation}
Therefore, the Green-Lagrange strain tensor is expressed as
\begin{equation}
\mathbf{E} = \left[ \varepsilon_{\alpha\beta}+\kappa_{\alpha\beta} \theta^{3} \right] \bar{\mathbf{g}}^{\alpha} \otimes \bar{\mathbf{g}}^{\beta} + \frac{1}{2} [\lambda_3^2 - 1] \bar{\mathbf{g}}^{3} \otimes \bar{\mathbf{g}}^{3},
\label{eq:GL_strian_1}
\end{equation}
with the mid-surface strain and curvature
\begin{equation}
 \varepsilon_{\alpha\beta} = \frac{1}{2} [a_{\alpha \beta} - \bar{a}_{\alpha \beta}]\quad \text{and} \quad \kappa_{\alpha\beta} =[- b_{\alpha \beta} + \bar {b}_{\alpha \beta}],
 \label{eq:GL_strian_2}
\end{equation}
corresponding to the stretching and bending strains, respectively.
The components of the Cauchy-Green and Green-Lagrange strain tensors are denoted as $C_{ij}$ and $E_{ij}$, respectively, in the following sections.

\subsection{Constitutive relations}
\label{sec:constitutive}
The Piola-Kirchhoff stress tensor is  a conjugate variable to the Green-Lagrange strain tensor. It is defined in terms of the covariant base vectors in the reference configuration as
\begin{equation}
    \mathbf S = S^{ij}\bar{\mathbf g}_i \otimes \bar{\mathbf g}_j.
\end{equation}
The components of the total differential of the Piola-Kirchhoff stress tensor, i.e. the elastic stiffness tensor, can be computed using the chain rule as
\begin{equation}
\mathrm{d} S^{ij} = \frac{\partial S^{ij}}{\partial E_{kl}} \mathrm{~d} E_{kl}=2 \frac{\partial S^{ij}}{\partial C_{kl}} \mathrm{~d} E_{kl} = \mathbb{C}^{ijkl} \mathrm{~d} E_{kl},
\label{eq:d_stress}
\end{equation}
where $\mathbb{C}^{ijkl}$ are the components of the fourth order elastic stiffness tensor $\boldsymbol{\mathbb{C}}$. 

We are concerned with incompressible hyperelastic solids in this work.
The incompressibility constraint is given by
\begin{equation}
    \text{det} (\mathbf F) = 1.
    \label{eq: J=1}
\end{equation}
The Piola-Kirchhoff stress tensor for incompressible hyperelastic solids~\cite{holzapfel2002nonlinear} is expressed as
\begin{equation}
    S^{ij} = 2\frac{\partial W }{\partial C_{ij}} - \tilde{p} {C}^{ij},
    \label{eq:stress_expression_0}
\end{equation}
where $C_{ij} $ are the covariant components of the Cauchy-Green strain tensor $\mathbf C$, while ${C}^{ij}$ are the contravariant components of ${\mathbf{C}}^{-1}$. 
For the present shell formulation, $C_{ij} = g_{ij}$ and ${C}^{ij} = {g}^{ij}$ as defined in Equations~\eqref{eq:Cauchy-Green_2} and~\eqref{eq:Cauchy-Green_inv}, respectively. 
The Lagrange multiplier $\tilde{p}$ enforcing incompressiblity is identified as the hydrostatic pressure inside the hyperelastic solid. 
\subsubsection{Plane stress condition for thin-shells}
{
The three-dimensional constitutive equations above can be simplified to two-dimensions by using the plane stress condition for thin shells $(S^{33} = 0)$ coupled with the incompressibility constraint \eqref{eq: J=1}.
Based on the Kirchhoff-Love assumption for thin shells, we neglect
the shear  strain components thereby writing
 the various components of Cauchy--Green deformation tensor and its inverse  as
\begin{subequations}
\begin{align}
C_{\alpha\beta} &= g_{\alpha\beta}, \quad C_{\alpha 3} = C_{3 \alpha} = 0 \quad \text{and} \quad C_{33} = g_{33} = \lambda_3^{2}, \\ 
C^{\alpha\beta} &= g^{\alpha\beta}, \quad {C}^{\alpha 3} = {C}^{3 \alpha} = 0 \quad \text{and} \quad {C}^{33} = {g}^{33} = \lambda_3^{-2} .
\end{align}
\end{subequations}
Using equation \eqref{eq:d_stress}, we write 
the in-plane stress components $S^{\alpha\beta}$ as
\begin{equation}
    S^{\alpha\beta} = 2\frac{\partial W }{\partial C_{\alpha\beta}} - \tilde{p} {C}^{\alpha\beta},
    \label{eq:stress_expression}
\end{equation}
and their differential as
\begin{equation}
\mathrm{d} S^{\alpha \beta} =2 \frac{\partial S^{\alpha \beta}}{\partial C_{\gamma\delta}} \mathrm{~d} E_{\gamma\delta} + 2 \frac{\partial S^{\alpha \beta}}{\partial C_{33}} \mathrm{~d} E_{33} = \mathbb{C}^{\alpha\beta\gamma\delta}  \mathrm{~d} E_{\gamma\delta} + \mathbb{C}^{\alpha\beta33}  \mathrm{~d} E_{33}.
\label{eq:d_stress_p}
\end{equation}
 For thin shells, the plane stress condition $S^{33} = 0$ results in
 \begin{equation}
     S^{33} = 2\frac{\partial W }{\partial C_{33}} - C^{33}\tilde{p} = 0, \quad \text{with} \quad C_{33} = \lambda_3^{2},\, C^{33} = \lambda_3^{-2}.
 \end{equation}
Thus $\tilde{p} $ can be explicitly determined as
\begin{equation}
    \tilde{p} = 2\frac{\partial W}{\partial C_{33}}{C}_{33} = 2\lambda_3^2\frac{\partial W}{\partial C_{33}}.
    \label{eq:lagrange_multiplier}
\end{equation}
Due to the plane-stress condition above, $\tilde p =\tilde{p}(C_{ij})$ thus can be considered as a function of the strain tensor $\mathbf E$ or the right Cauchy--Green deformation tensor $\mathbf{C}$. Thus the total derivative of the Piola--Kirchhoff stress tensor is now written as
\begin{equation}
\frac{\mathrm d }{\mathrm d \mathbf E}\big(\mathbf S\left(\mathbf E, \tilde{p}(\mathbf E)\right)\big) = \frac{\partial \mathbf S}{\partial \mathbf E} + \frac{\partial \mathbf S}{\partial \tilde p} \frac{\partial \tilde p}{\partial \mathbf E} =  2\frac{\partial \mathbf S}{\partial \mathbf C} + 2\frac{\partial \mathbf S}{\partial \tilde p} \frac{\partial \tilde p}{\partial \mathbf C}.
\end{equation}
Upon comparison with Equations~\eqref{eq:d_stress} and~\eqref{eq:stress_expression} and using \eqref{eq:lagrange_multiplier}, one can write the explicit expression for the in-plane components of the fourth-order tensor $\boldsymbol{\mathbb{C}}$ as (see \cite{kiendl2015isogeometric} for detailed derivations)
\begin{equation}
\mathbb{C}^{\alpha\beta\gamma\delta}  = 4\frac{\partial^2 W }{\partial C_{\alpha \beta} \partial C_{\gamma\delta}} 
- 2\frac{\partial \tilde{p}}{\partial C_{\gamma\delta}}{C}^{\alpha\beta} - 2\frac{\partial \tilde{p}}{\partial C_{\alpha\beta}}{C}^{\gamma\delta}  - \tilde{p} [ {C}^{\alpha\beta}  {C}^{\gamma\delta} -  {C}^{\alpha\gamma} {C}^{\beta\delta} - {C}^{\alpha\delta}  {C}^{\beta\gamma}],
\end{equation}
}
where, the partial derivative of the Lagrange multiplier can be calculated from Equation~\eqref{eq:lagrange_multiplier} as
\begin{equation}
\frac{\partial \tilde{p}}{\partial C_{\alpha\beta}} = 2\frac{\partial^2 W}{\partial C_{33} \partial C_{\alpha\beta}}\lambda_3^2.
\end{equation}
Furthermore, the plane stress condition requires the incremental stress in the thickness direction to vanish, that is
\begin{equation}
\mathrm d S^{33} = \mathbb{C}^{33\alpha\beta} \mathrm d  E_{\alpha\beta} + \mathbb{C}^{3333}\mathrm d  E_{33} = 0. 
\end{equation}
This relation can be used to explicitly calculate the differential strain component $\mathrm d E_{33} = -\mathbb{C}^{33\alpha\beta} \mathrm d E_{\alpha\beta} / \mathbb{C}^{3333} $ and upon substituting it into Equation~\eqref{eq:d_stress_p}, the in-plane tangent tensor $\boldsymbol{\mathbb{C}}$ is modified as
\begin{align}
\hat{\mathbb C}^{\alpha\beta\gamma\delta} = \mathbb{C}^{\alpha\beta\gamma\delta}  - \frac{ \mathbb C^{\alpha\beta33}  \mathbb C^{33\gamma\delta}}{ \mathbb C^{3333}},
\label{eq:mod_elastic_tensor}
\end{align}
where
\begin{align}
\mathbb{C}^{\alpha \beta 33} &= 4\frac{\partial^2 W }{\partial C_{\alpha \beta} \partial C_{33}} 
- 2\frac{\partial \tilde{p}}{\partial C_{33}}{C}^{\alpha\beta} - 2\frac{\partial \tilde{p}}{\partial C_{\alpha\beta}}{C}^{33}  - \tilde{p} [ {C}^{\alpha\beta}  {C}^{33}].
\label{eq:Cab33}
\end{align}
The partial derivative of the Lagrange multiplier with respect to $C_{33}$ is given as
\begin{equation}
    \frac{\partial \tilde{p}}{\partial C_{33}} = 2\frac{\partial^2 W}{\partial C_{33}^2}\lambda_3^2 + 2 \frac{\partial W}{\partial C_{33}}.
\end{equation}
Thus Equation~\eqref{eq:Cab33} can be computed explicitly as 
\begin{align}
\mathbb{C}^{\alpha \beta 33} = -{C}^{\alpha \beta} \left[ 6 \frac{\partial W}{\partial C_{33}} +4 \frac{\partial^{2} W}{\partial C_{33}^{2}} \lambda_3^2 \right].
\label{eq: cab33}
\end{align}
Similarly, the other two terms in Equation~\eqref{eq:mod_elastic_tensor} are
\begin{align}
\mathbb{C}^{33\gamma \delta} = -{C}^{\gamma \delta} \left[ 6 \frac{\partial W}{\partial C_{33}} +4 \frac{\partial^{2} W}{\partial C_{33}^{2}} \lambda_3^2 \right], \quad
\mathbb{C}^{3333} = -\lambda_3^{-2} \left[ 6 \frac{\partial W}{\partial C_{33}} +4 \frac{\partial^{2} W}{\partial C_{33}^{2}} \lambda_3^2 \right],
\label{eq: c33gd c3333}
\end{align}
and therefore a closed form expression for $\hat{\mathbb{C}}^{\alpha \beta \gamma \delta}$ is obtained by substituting \eqref{eq: cab33} and \eqref{eq: c33gd c3333} in Equation \eqref{eq:mod_elastic_tensor}.

\subsubsection{Incompressible Mooney-Rivlin material}
Due to incompressibility, volume of the shell remains unchanged, that is
\begin{equation}
\int_{\bar\Gamma} \bar{h} \, \mathrm d \bar{\Gamma} = \int_{\Gamma} \lambda_3 \bar{h} \, \mathrm d {\Gamma},
\end{equation}
where $\mathrm d \Gamma = \mathcal J \mathrm d \bar{\Gamma}$ and $\mathcal J$ is the in-plane Jacobian, which is expressed as
\begin{equation}
\mathcal J  =\frac{|\mathbf a_1 \times \mathbf a_2|}{|\bar{\mathbf a}_1 \times \bar{\mathbf a}_2|} = \frac{J}{\bar J}.
\label{eq:cal_J}
\end{equation}
Thus the thickness stretch is the inverse of in-plane Jacobian, that is
\begin{equation}
\lambda_3 = \mathcal J^{-1} =\frac{|\bar{\mathbf{a}}_1 \times \bar{\mathbf{a}}_2|}{|\mathbf{a}_1 \times \mathbf{a}_2|}.
\end{equation}
The elastic solid is modelled using the Mooney--Rivlin hyperelastic constitutive law.
This strain energy density per unit undeformed volume is expressed as
\begin{equation}
    W(\mathbf{C}) = c_1[I_1 - 3] + c_2[I_2 - 3],
\end{equation}
where $I_1$ and $I_2$ are the first and second invariants of $\mathbf{C}$, defined by
\begin{equation}
I_1 = \text{tr}(\mathbf C)= \lambda_1^2 + \lambda_2^2 + \lambda_3^2\quad \text{and} \quad I_2 = \frac{1}{2}\left[ \text{tr}(\mathbf C)^2 -\text{tr}(\mathbf C^2) \right]= \lambda_1^{-2} + \lambda_2^{-2} + \lambda_3^{-2}.
\end{equation}
If $c_2 = 0$, the model reduces to neo-Hookean material. We explicitly calculate the derivatives of the strain energy density for future use as
\begin{subequations}
\begin{align}
    \frac{\partial W}{\partial C_{\alpha\beta}} = c_1  \bar{g}^{\alpha\beta} + c_2[C_{\gamma\delta}\bar{g}^{\gamma\delta} \bar{g}^{\alpha\beta} - \bar{g}^{\alpha \gamma}C_{\gamma\delta}\bar{g}^{\delta\beta}],
    \quad 
    \frac{\partial W}{\partial C_{33}} = c_1 + c_2[C_{\gamma\delta}\bar{g}^{\gamma\delta} - \lambda_3^2], \\
    \frac{\partial^2 W }{\partial C_{\alpha\beta} \partial C_{\gamma\delta}} = c_2\bar{g}^{\gamma\delta} \bar{g}^{\alpha\beta} - \frac{c_2}{2} [\bar{g}^{\alpha\gamma} \bar{g}^{\beta\delta} + \bar{g}^{\alpha\delta} \bar{g}^{\beta\gamma}],
    \quad 
    \frac{\partial^2 W }{\partial C_{\alpha\beta} \partial C_{33}} = c_2 \bar{g}^{\alpha\beta},
    \quad
    \frac{\partial^{2} W}{\partial C_{33}^{2}} = 0.
\end{align}
\end{subequations}

\subsection{Stress resultants for thin shells}
The thin shell formulation considers a three dimensional solid as a two dimensional surface with a thickness. Thus, the internal forces are the stress resultants integrated through the thickness. These are decomposed into the normal force $\hat{\mathbf n}$ and the bending moment $\hat{\mathbf m} $, with their components calculated from the Piola-Kirchhoff stress as
\begin{equation}
\begin{aligned}
n^{\alpha \beta} &=\int_{-\frac{\bar h}{2}}^{\frac{\bar h}{2}} S^{\alpha \beta} J_c \mathrm{~d} \theta^{3}, \\
m^{\alpha \beta} &=\int_{-\frac{\bar h}{2}}^{\frac{\bar h}{2}} S^{\alpha \beta} \theta^{3} J_c \mathrm{~d} \theta^{3},
\end{aligned}
\end{equation}
where 
\begin{equation}
J_c = \frac{|[\bar{\mathbf{g}}_1 \times \bar{\mathbf{g}}_2] \cdot \bar {\mathbf{g}}_3|}{|[\bar {\mathbf{a}}_1 \times \bar {\mathbf{a}}_2] \cdot \bar {\mathbf{a}}_3|}
\end{equation}
denotes the change of the in-plane Jacobian along the thickness direction. Based on~\eqref{eq:GL_strian_1} and~\eqref{eq:GL_strian_2} for $E_{\alpha \beta}$, their total differentials are calculated as
\begin{subequations}
\begin{align}
\mathrm{d} n^{\alpha \beta} &=\left[\int_{-\frac{\bar h}{2}}^{\frac{\bar h}{2}} \hat{\mathbb{C}}^{\alpha \beta \gamma \delta} J_c \mathrm{~d} \theta^{3}\right] \mathrm {~d} \varepsilon_{\gamma \delta}+\left[\int_{-\frac{\bar h}{2}}^{\frac{\bar h}{2}} \hat{\mathbb{C}}^{\alpha \beta \gamma \delta} \theta^{3} J_c \mathrm{~d} \theta^{3}\right] \mathrm{~d} \kappa_{\gamma \delta} ,\\
\mathrm{d} m^{\alpha \beta} &=\left[\int_{-\frac{\bar h}{2}}^{\frac{\bar h}{2}} \hat{\mathbb{C}}^{\alpha \beta \gamma \delta} \theta^{3} J_c \mathrm{~d} \theta^{3}\right] \mathrm{~d} \varepsilon_{\gamma \delta}+\left[\int_{-\frac{\bar h}{2}}^{\frac{\bar h}{2}} \hat{\mathbb{C}}^{\alpha \beta \gamma \delta}\left[\theta^{3}\right]^{2} J_c \mathrm{~d} \theta^{3}\right] \mathrm {~d} \kappa_{\gamma \delta} .
\end{align}
\end{subequations}
\subsection{Virtual work principle}
\label{sec:virtual_work}
The virtual work principle is applied in order to obtain the governing equations for the nonlinear thin shells theory
\begin{equation}
\delta \mathcal W_\text{int} - \delta \mathcal W_\text{ext} = 0,
\label{eq:int_vw_varition}
\end{equation}
The internal virtual work is given by
\begin{equation}
\delta \mathcal W_\text{int} = \int_{\bar\Gamma} \left[n^{\alpha\beta} \delta \varepsilon_{\alpha\beta} +  m^{\alpha\beta} \delta \kappa_{\alpha\beta}\right]\mathrm {~d} \bar\Gamma,
\label{eq:ext_vw_varition}
\end{equation}
and the external virtual work by
\begin{equation}
    \delta \mathcal W_\text{ext} =\hat{h}\int_\Gamma  b_i \delta u_i \,\mathrm d \Gamma + \hat{h} \int_{S_t} \tau_i \delta u_i \, \mathrm d S_t, 
\end{equation}
with
\begin{equation}
    \hat{h} = \int_{-\frac{h}{2}}^{\frac{h}{2}} J_c \mathrm{~d} \theta^3.
\end{equation}
$b_i$ denotes the components of a body force. The components $\tau_i$ are prescribed tractions at the boundary $S_t\in \partial \Gamma$ and vanish for enclosed geometries. 
For convenience, the internal virtual work is integrated in the reference configuration, while the external virtual work is calculated in the deformed configuration.
\section{Numerical Implementation}
\label{sec:implementation}
\subsection{Catmull-Clark subdivision surfaces}
The fundamental idea of subdivision surfaces is to generate a smooth surface by repeatedly refining a coarse control grid using a subdivision scheme {that generalizes bi-cubic uniform B-spline knot insertion. For arbitrary initial meshes, this scheme generates limit surfaces that are $C^2$ continuous everywhere except at extraordinary vertices where they guarantee a $C^1$ limiting surface.}
For regular vertices, the limit surface generated by the Catmull-Clark subdivision scheme is identical to a bi-cubic B-spline surface. 
\autoref{fig:subd_surface} shows a patch of a subdivision surface with its control grid for regular vertices. The grid divides the parametric domain of the surface patch into nine elements, see \autoref{fig:CC_element_patch}. The surface point $\bar{\mathbf x}$ in the central element can be interpolated using a tensor product of two cubic B-splines with 16 control points as
\begin{equation}
\bar{\mathbf x}(\theta^1, \theta^2) = \sum_{a=0}^{15} \tilde{N}_a(\theta^1,\theta^2)\mathbf P_{a},
\end{equation}
where $\mathbf P_{a}$ is the $a^{\text{th}}$ control point and $\tilde{N}_a$ denotes the corresponding local base function for the element which is defined by
\begin{equation*}
    \tilde{N}_a(\theta^1,\theta^2) = N_{a\%4}(\theta^1)N_{\lfloor{a/4}\rfloor}(\theta^2),
\end{equation*}
here $\lfloor \bullet \rfloor$ is the modulus operator and $\%$ denotes the remainder operator which gives the remainder of the integer division.
\begin{figure}[h]
\centering
\begin{subfigure}[b]{0.5\linewidth}
	\centering
 	 \includegraphics[width=\linewidth]{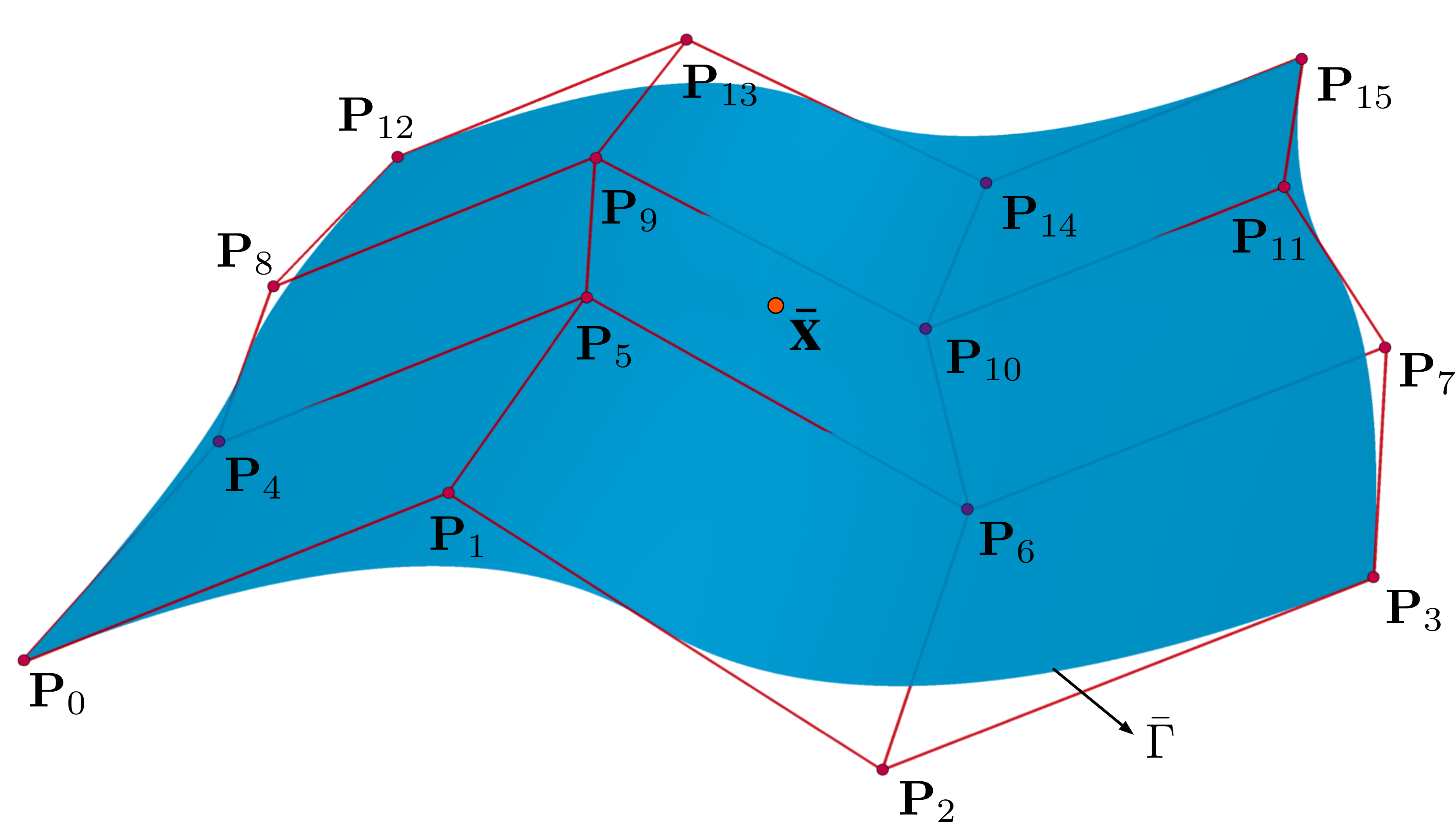}
	\caption{}
	\label{fig:subd_surface}
\end{subfigure}
\begin{subfigure}[b]{0.35\linewidth}
	\centering
 	 \includegraphics[width=\linewidth]{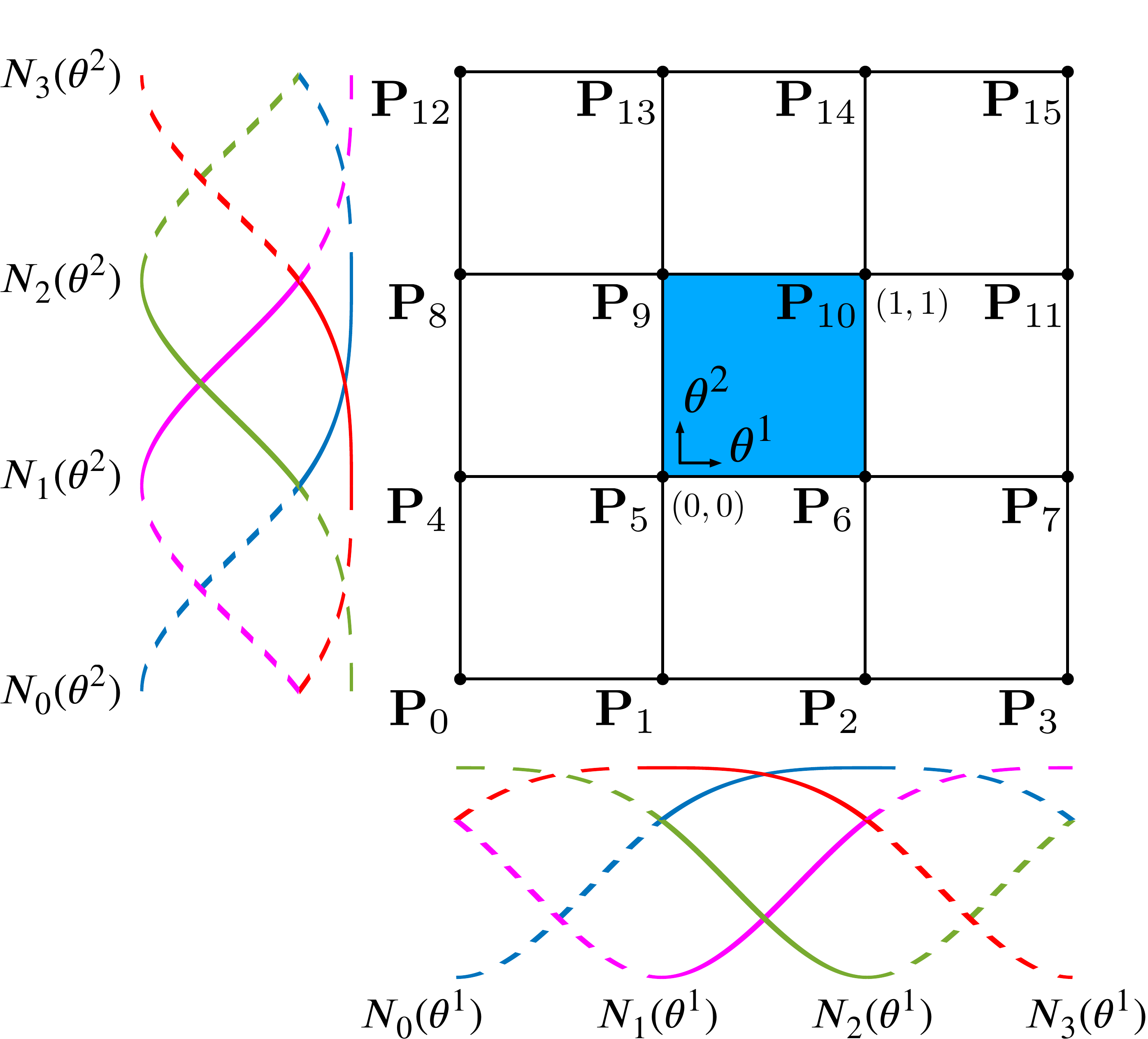}
	\caption{}
	\label{fig:CC_element_patch}
\end{subfigure}
\caption{(a) An example patch of a Catmull-Clark subdivision surface. (b) The parametric domain of the patch and the corresponding bases functions.}
\label{fig:subd_surfaces}
\end{figure}
As a result of the tensor-product nature, each vertex in a subdivision surface control grid is connected with only four elements. One defines the number of elements connected as the `valence' of the vertex. A regular vertex in a Catmull-Clark surface control grid has a valence of 4. However, contrary to NURBS surfaces, Catmull-Clark subdivision surfaces can handle irregular cases where the valence of a vertex is not equal to 4. Thus it  allows to handle complex geometry with arbitrary topology. These vertices are known as `extraordinary vertices' and require a different algorithm for the evaluation of the limiting surface, introduced in~\cite{stam1998exact, liu2020assessment}.  Catmull-Clark subdivision surfaces display $C^1$ continuity at extraordinary vertices~\cite{peters1998analysis}, otherwise they possess $C^2$ continuity. Therefore, Catmull-Clark subdivision surfaces provide an adequate $C^1$ continuous discretisation to satisfy the requirement of the Galerkin formulation of Kirchhoff-Love shells, where its test and trial functions must be in the Hilbert space $H^2(\Omega)$~\cite{cirakortiz2000}.
\subsection{Discretisation and linearisation}
The displacement of the mid-surface is discretised using Catmull-Clark subdivision bases as shape functions:
\begin{equation}
\mathbf u = \sum_{A = 1} ^{n_b} N^{A} \mathbf u^{A},
\end{equation}
where $n_b$ is the total number of basis functions and is equal to the total number of the control points, $N^{A}$ is the basis function corresponding to the $A^{\text{th}}$ control point. We note here that $A$ is the global index of the control point. $\mathbf u^{A}$ denotes the $A^{\text{th}}$ nodal displacement vector with three components corresponding to three Cartesian coordinates denoted as $u^{A}_i$, leading to the total number of degrees of freedom to interpolate $\mathbf u$ equal to $3n_b$. The variation of $\mathbf u$ with respect to its components is expressed as
\begin{equation}
\frac{\partial \mathbf u}{\partial u^A_i} = N^{A} \mathbf c_i,
\end{equation}
where $\mathbf c_i$ are the orthonormal basis vectors.
Index $r$ varies from $1$ to $3n_b$ referring the degree of freedom and $r = 3[A-1] + i$, such that $u_r$ denotes the $i^\text{th}$ component of $\mathbf u^A$. Thus, the variation of $\mathbf u$ with respect to the $r^\text{th}$ degree of freedom reads
\begin{equation}
    \delta_{r}\mathbf u = N^{A} \mathbf c_{i}.
\end{equation}
The variation of the strain components $\varepsilon_{\alpha\beta}$ and $\kappa_{\alpha\beta}$ are expressed as
{  
\begin{equation}
\begin{aligned}
 \delta_{r}\varepsilon_{\alpha\beta} &= \frac{1}{2} [ \delta_{r}{\mathbf a_{\alpha}} \cdot \mathbf a_{\beta} +   \delta_{r}{\mathbf{a}_{\beta}} \cdot \mathbf a_{\alpha}],\\
 \delta_{r}\kappa_{\alpha\beta} 
 &=  -  \lambda_3 [\delta_{r} \mathbf {a_{\alpha,\beta}}\cdot \mathbf a_3 + \mathbf a_{\alpha,\beta}\cdot  \delta_{r}\mathbf a_{3}] - \delta_r \lambda_3 [\mathbf{a}_{\alpha,\beta} \cdot \mathbf a_3].\\
\end{aligned}
\end{equation}
}
The variation of $\mathbf a_{\alpha}$ and $\mathbf a_{\alpha,\beta}$ are easily expressed using the first and second derivatives of the basis functions as
\begin{equation}
 \delta_{r}{\mathbf a_{\alpha}} = N^{A}_{,\alpha} \mathbf c_i, \quad
 \delta_{r}{\mathbf a_{\alpha,\beta}} = N^A_{,\alpha\beta} \mathbf c_i.
\end{equation}
The variations of the normal vector $ \delta_{r} \mathbf a_3$ {  and the thickness stretch $\delta_r \lambda_3$ } can be found in the appendix.
The variations of the internal and external virtual work with respect to the $r^\text{th}$ degree of freedom are then defined by
\begin{equation}
\begin{aligned}
\delta_r \mathcal W_{\text{int}} &= \int_{\bar\Gamma} \left[n^{\alpha\beta} \delta_{r}\varepsilon_{\alpha\beta} + m^{\alpha\beta} \delta_{r}\kappa_{\alpha\beta} \right] \mathrm{~d} \bar\Gamma,\\
\delta_r \mathcal W_{\text{ext}} &=\hat{h}\int_\Gamma  b_i \delta_r u_i \,\mathrm d \Gamma,
\end{aligned}
\end{equation}
where the boundary traction has been ignored for simplicity. The global residual vector $\ary R$ is defined by
\begin{equation}
\ary R = \ary F^{\text{int}} - \ary F^{\text{ext}},
\end{equation}
where $\ary F^{\text{int}}$ and $\ary F^{\text{ext}}$ are two global vectors with $3n_b$ components each, and where the $r^\text{th}$ components are $\delta_r \mathcal W_{\text{int}}$ and $\delta_r \mathcal W_{\text{ext}}$, respectively. The system is in equilibrium when the residual vector $\ary R = \ary 0$. The global tangential stiffness matrix $\ary{K}$ of $3n_b \times 3n_b$ has entries in $r^{\text{th}}$ row and $s^{\text{th}}$ column given by
\begin{equation}
 K^{rs} = \delta_s\delta_r \mathcal W_{\text{int}} - \delta_s\delta_r \mathcal W_{\text{ext}},
\end{equation}
where
\begin{equation}
\begin{aligned}
\delta_s\delta_r \mathcal W_{\text{int}} &= \int_{\bar\Gamma} \left[\delta_{s} n^{\alpha\beta} \, \delta_{r}\varepsilon_{\alpha\beta} + n^{\alpha\beta} \, \delta_{s}\delta_{r}\varepsilon_{\alpha\beta} + \delta_{s} m^{\alpha\beta} \, \delta_{r}\kappa_{\alpha\beta} + m^{\alpha\beta} \, \delta_{s}\delta_{r}\kappa_{\alpha\beta} \right] \mathrm{~d} \bar\Gamma \\
\delta_s\delta_r \mathcal W_{\text{ext}} &= \hat{h}\int_\Gamma \delta_s  b_i  \,\delta_r u_i \,\mathrm d \Gamma,
\end{aligned}
\end{equation}
and
\begin{equation}
\begin{aligned}
\delta_s n^{\alpha \beta} &=\left[\int_{-\bar{h} / 2}^{\bar{h}  / 2} \hat{\mathbb{C}}^{\alpha \beta \gamma \delta} J_c \mathrm{~d} \theta^{3}\right] \delta_s \varepsilon_{\gamma \delta}+\left[\int_{-\bar{h} / 2}^{\bar{h} / 2} \hat{\mathbb{C}}^{\alpha \beta \gamma \delta} \theta^{3} J_c \mathrm{~d} \theta^{3}\right] \delta_s \kappa_{\gamma \delta} \\
\delta_s  m^{\alpha \beta} &=\left[\int_{-\bar{h} / 2}^{\bar{h} / 2} \hat{\mathbb{C}}^{\alpha \beta \gamma \delta} \theta^{3} J_c \mathrm{~d} \theta^{3}\right] \delta_s \varepsilon_{\gamma \delta}+\left[\int_{-\bar{h} / 2}^{\bar{h} / 2} \hat{\mathbb{C}}^{\alpha \beta \gamma \delta}\left[\theta^{3}\right]^{2} J_c \mathrm{~d} \theta^{3}\right]\delta_s \kappa_{\gamma \delta},
\end{aligned}
\end{equation}
and
{ 
\begin{equation}
\begin{aligned}
\delta_{s}\delta_{r}\varepsilon_{\alpha\beta} =& \frac{1}{2} [ \delta_{r}{\mathbf a_{\alpha}} \cdot \delta_{s} \mathbf a_{\beta} +   \delta_{r}{\mathbf{a}_{\beta}} \cdot \delta_{s} \mathbf a_{\alpha}], \\
\delta_{s} \delta_{r}\kappa_{\alpha\beta} = & -  \delta_{s}\lambda_3 [\delta_{r} \mathbf {a_{\alpha,\beta}}\cdot \mathbf a_3 + \mathbf a_{\alpha,\beta}\cdot  \delta_{r}\mathbf a_{3}]  \\
&- \lambda_3 [
\delta_{r} \mathbf {a_{\alpha,\beta}}\cdot \delta_{s} \mathbf a_3 + \delta_{s} \mathbf a_{\alpha,\beta}\cdot  \delta_{r}\mathbf a_{3} + \mathbf a_{\alpha,\beta}\cdot  \delta_s\delta_{r}\mathbf a_{3}] \\
& - \delta_s \delta_r \lambda_3 [\mathbf{a}_{\alpha,\beta} \cdot \mathbf a_3] - \delta_r \lambda_3 [\delta_s\mathbf{a}_{\alpha,\beta} \cdot \mathbf a_3 + \mathbf{a}_{\alpha,\beta} \cdot \delta_s\mathbf a_3].\\
\end{aligned}
\end{equation}
}
The second variation of the normal vector $ \delta_{s}\delta_{r} \mathbf a_3$ and { $\delta_s \delta_r \lambda_3 $} can be found in the Appendix~\ref{sec:Appendix_1}.
Here, it is important to note that for inflated thin shells, the external load $\ary F_\text{ext}$ is a follower load and is a function of the displacements. 
Hence the load increment must consider the deformation of the geometry and the tangential term $\delta_s\delta_r \mathcal W_{\text{ext}}$ should be included in the tangential stiffness matrix. 
For example, if the external load is a uniform pressure inside an enclosed hyperelastic thin shell, the first variation of the external virtual work with respect to the degrees of freedom is given by
\begin{equation}
\delta_r \mathcal W_{\text{ext}} = \int_{\Gamma} p [\mathbf a_3 \cdot \delta_r \mathbf u] \mathrm{~d} \Gamma,
\end{equation}
where $\mathrm d \Gamma = \mathcal J \mathrm d \bar{\Gamma}$. Substituting the expression~\eqref{eq:cal_J} into the above equation allows the second variation of the external virtual work with respect to the degrees of freedom to be derived as
\begin{equation}
\delta_s\delta_r \mathcal W_{\text{ext}} = \int_{\bar\Gamma} p \left[[\delta_{s}\mathbf a_1 \times \mathbf a_2 + \mathbf a_1 \times \delta_{s}\mathbf a_2] \cdot \delta_r \mathbf u\right] \bar{J}^{-1} \mathrm{~d} \bar\Gamma.
\end{equation}
 Here, the pressure $p$ is considered as an external load, and the global vector $\ary F_{\text{ext}}$  must be regenerated for each load step.

\section{Large deformation and associated instabilities}
\label{sec:investigations}
\subsection{Nonlinear algorithm}
The inflated thin structures considered in the present work are incompressible hyperelastic shells under uniform inflating pressure. 
The inflation of such soft and thin structures causes large deformations and a highly nonlinear mechanical response. 
The hyperelastic shell formulation involves both kinematic and constitutive nonlinearities.
\begin{figure}[]
\centering
  \includegraphics[width=0.8\linewidth]{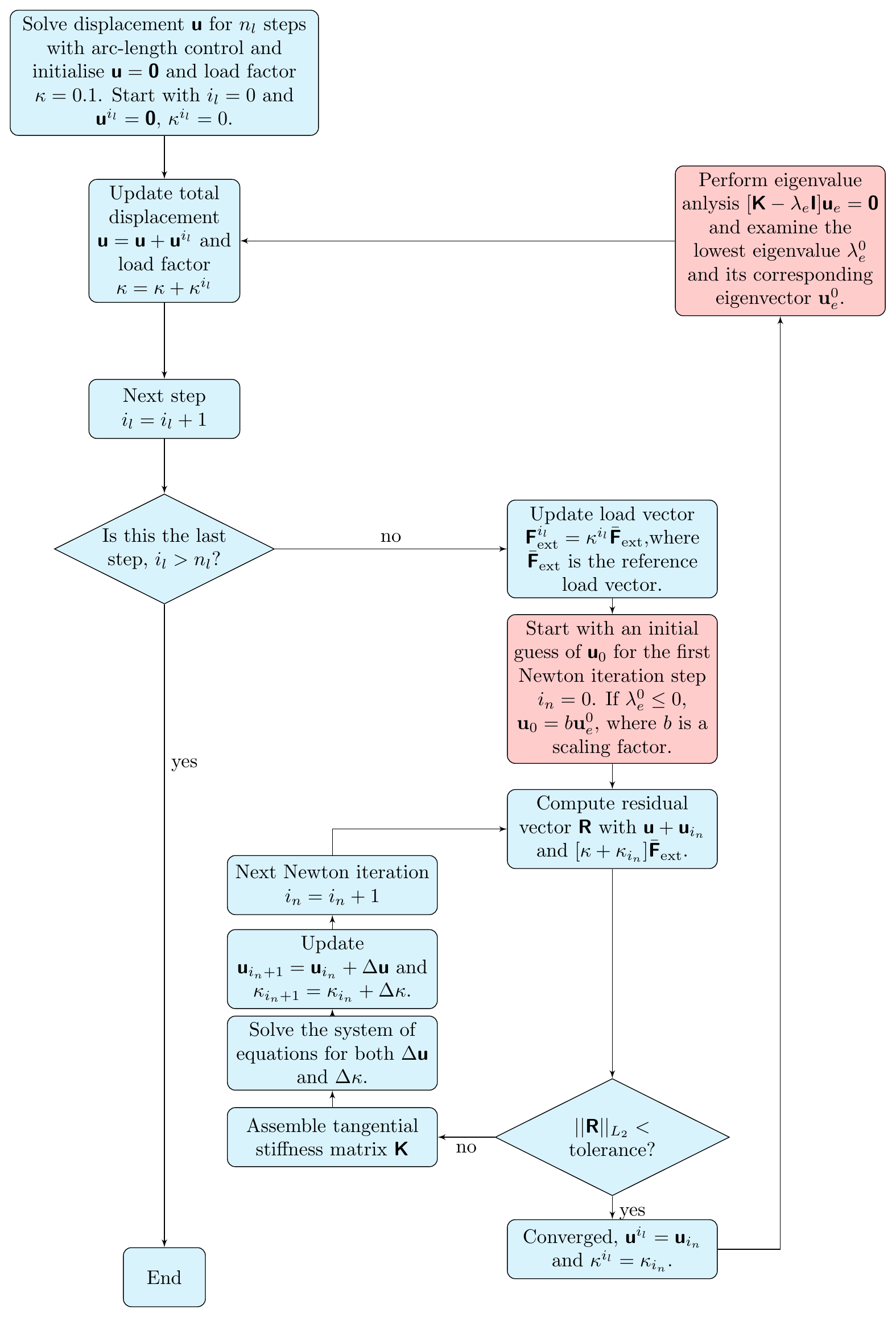}
\caption{Algorithm for solving the nonlinear hyperelastic thin shell formulation.}
\label{fig:nonlinear_algorithm}
\end{figure}
\autoref{fig:nonlinear_algorithm} shows a flowchart illustrating the adopted algorithm for such nonlinear problems. Since the mechanical response of the hyperelastic shell structures is nonlinear and can be quite complex, the problem is solved incrementally and arc-length method~\cite{kadapa2021simple} is adopted in the present work. For each increment, the Newton-Raphson method is employed to solve the nonlinear system of equations. It starts with an initial guess of the global nodal displacement vector $\ary u_0$. The linearised system of equations are formulated to solve the global vector of the incremental nodal displacements $\Delta\ary u$ and the incremental load factor $\Delta \kappa$ (see~\autoref{fig:nonlinear_algorithm}).
For each Newton iterate step, the displacement is updated until the $L_2$ norm of the residual vector is less than a tolerance which is set as $1\%$ of the initial norm of the residual vector. 

\subsection{Limit point instability}

\begin{figure}[]
\centering
\begin{subfigure}[b]{0.4\linewidth}
	\centering
 	 \includegraphics[width=\linewidth]{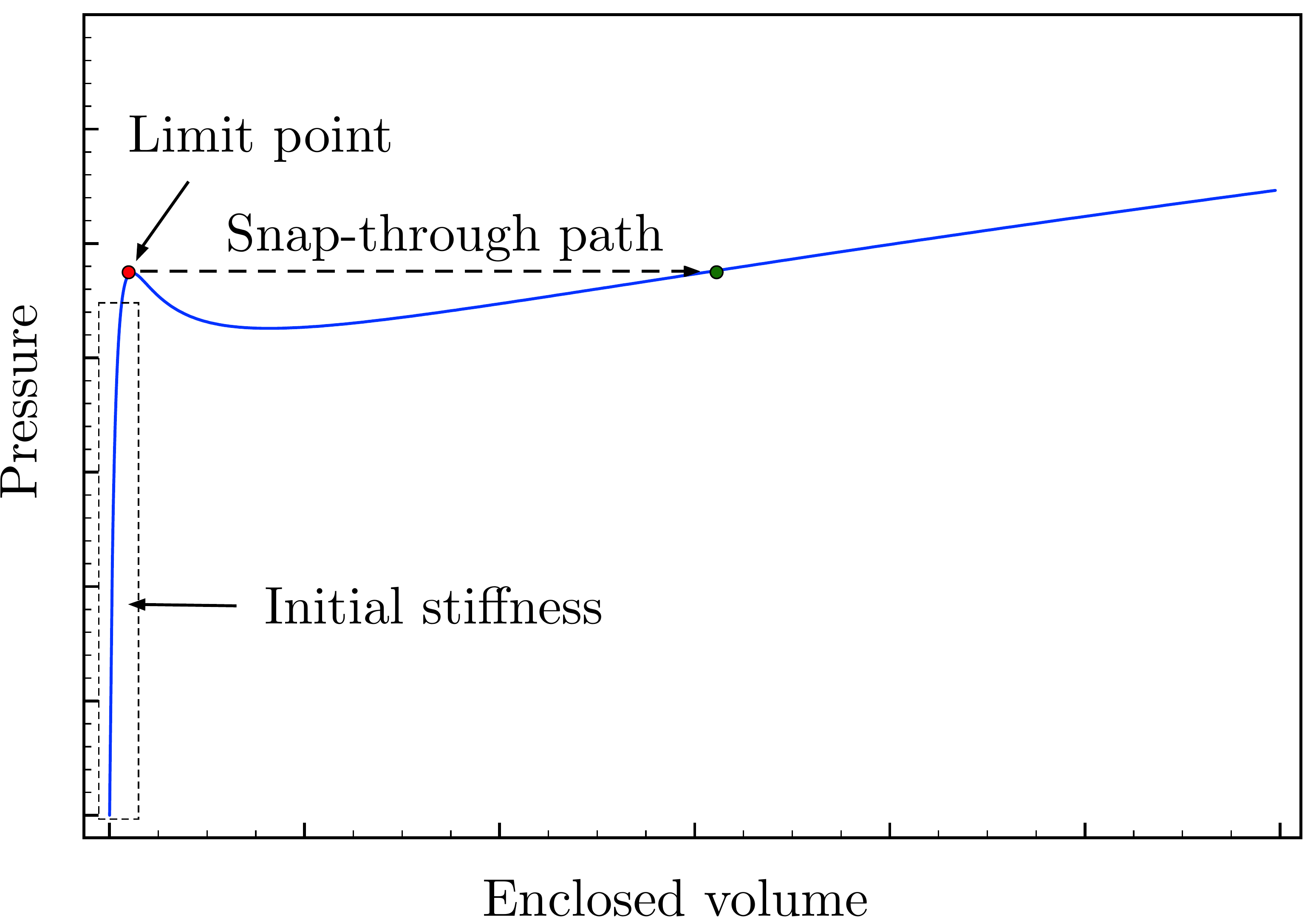}
	\caption{}
	\label{fig:snap_through}
\end{subfigure}
\begin{subfigure}[b]{0.4\linewidth}
	\centering
 	 \includegraphics[width=\linewidth]{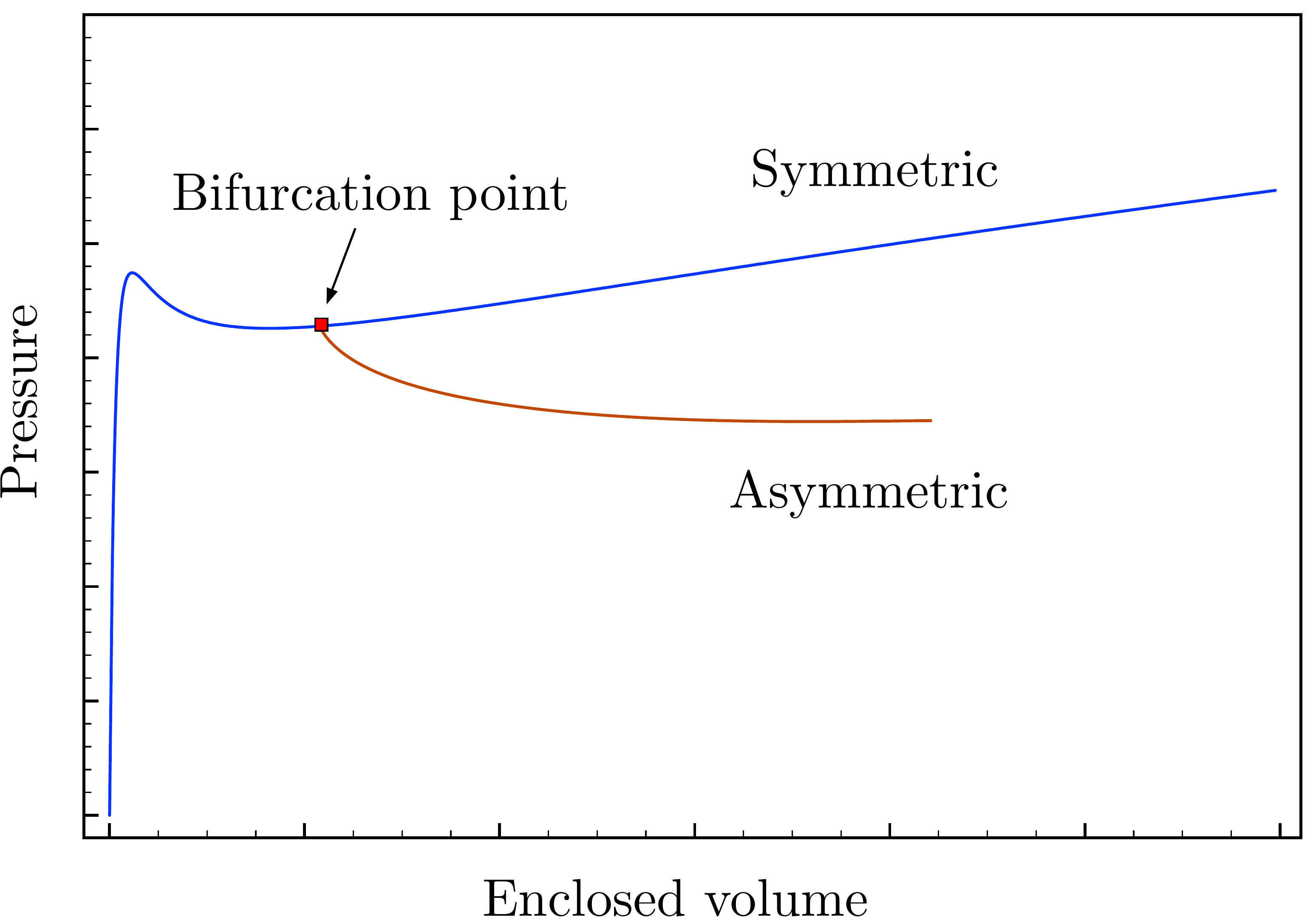}
	\caption{}
	\label{fig:bifurcation_sketch}
\end{subfigure}
\caption{Representative volume-pressure plots for hyperelastic thin shells illustrating two types of instabilities: (a) Snap-through and (b) Bifurcation.}
\label{fig:instabilities}
\end{figure}
The inflation of a hyperelastic thin shell will cause instability due to its highly nonlinear mechanical response. The most well-known instability phenomenon for inflated membranes/shells is snap-through buckling. \autoref{fig:snap_through} shows a plot of the enclosed volume against the internal pressure of an ideal hyperelastic shell. The structure has large initial stiffness when the internal pressure is relatively low and it is difficult to be inflated. With the increase of internal pressure, the enclosed volume changes accordingly, and the structure gradually loses its stiffness. When the pressure approaches the limit point, the structure dramatically loses its stiffness and undergoes very large inflation and jumps to a new equilibrium state. For a thin shell with imperfections, snap-through may occur prior to reaching the limit point. To investigate this type of instability, the present work uses a path-following and branching method to compute accurate assessment plot of pressure against volume for structures with post-buckling behaviour to capture the non-uniqueness of solution.

\subsection{Bifurcation from the principal solution}
\label{sec:eigenvalue}
For a closed hyperelastic shell with initial geometric symmetry, the inflated structure will remain symmetric until the internal volume exceeds a critical value, whereafter the structure may lose its symmetry under a small perturbation. This is due to the non-uniqueness of the solution of the nonlinear problem and the total energy for the non-symmetrical state being lower than that of the symmetric state. The point where the solution curve
differs is called the bifurcation point and is shown in \autoref{fig:bifurcation_sketch}. The bifurcation phenomenon is widely observed in mechanical experiments due to manufacturing imperfections and perturbations. However, numerically simulating bifurcation instability is difficult because it is a sudden change in a highly nonlinear problem. Moreover, a slight change of the geometry or material property may significantly affect the critical point at which the bifurcation occurs.
The proposed method to determine bifurcation assesses the stiffness matrix of the structure after solving each load increment. An accompanying eigenvalue analysis is performed by solving 
\begin{equation}
   [\ary{K}  -\lambda_e\ary{I}]\ary{u}_e = \ary{0},
\end{equation}
where $\lambda_e$ is the eigenvalue, $\ary I$ denotes an identity matrix, and the corresponding eigenvector is denoted by $\ary{u}_e$.
The dimensions of the stiffness matrix are $3n_b \times 3n_b$, thus a total of $3n_b$ eigenvalues can be computed. However, only the first few eigenvalues are critical as they are dominant. After eliminating the rigid body motion and rotation, if the smallest eigenvalue of the stiffness matrix is close to zero, the solution of the next load increment may switch into another solution branch. In the present work, the corresponding eigenvector is used as the initial guess for the next load step in order to perturb the structure and induce the bifurcation. It is noteworthy that the eigenvector is a normalised vector, which needs to be scaled in order to serve as a displacement perturbation. The choice of the scaling factor is vital: if the factor is too small, the solution will still follow the original path, and if the factor is too large, the system of equations will be difficult to solve. This method is embedded in our nonlinear algorithm and coloured in red in~\autoref{fig:nonlinear_algorithm}.

\section{Numerical examples}
In this section, we present four numerical examples to illustrate our method which has been implemented using the finite element library deal.II~\cite{Arndt2022-dealii94,Arndt2021-dealii-generic}.
We first solve the benchmark problem of the inflation of a Mooney--Rivlin circular plate to validate our method and demonstrate its accuracy.
Inflation of a spherical shell is computed for both neo-Hookean and Mooney--Rivlin models and validated with analytical solutions. 
This problem also demonstrates the ability of our implementation of the arc-length method to follow a nonliner path and capture the limit point instability.
Next, the inflation of a Mooney--Rivlin toroidal shell is computed. 
In this example, in addition to the limit point instability, we illustrate the computation of a bifurcation point and the post-bifurcation response of the toroidal shell.
Finally, we demonstrate the ability of our method to model complex deformation of an arbitrary geometry by computing the inflation of an airbag modelled using a Saint Venant--Kirchhoff constitutive law.
We are able to calculate large deformation of this shell structure along with complex folds and wrinkles.
All computations are performed using dimensionless quantities.

\subsection{Circular plate inflated with uniform pressure}
The first example considered is a circular inflated hyperelastic plate. The plate is simply supported and the pressure is considered to be a uniform load which is always perpendicular to the deformed mid-surface of the plate. Thus the pressure is a follower load and it is a function of the deformation. This problem was first analysed in~\cite{oden1970analysis} and it has become a benchmark problem for incompressible hyperelastic shells~\cite{hughes1983nonlinear,Cirak:2001aa,nama2020nonlinear}. \autoref{fig:plate_mesh} shows the initial control grid with $90$ elements and the limit surface of the circular plate. 
{
The radius of the plate is taken as $7.5$ while the thickness is $0.5$.
The material parameters for the Mooney--Rivlin model are $c_1 = 80 \mu, c_2 = 20\mu$, where $\mu$ is the shear modulus.
We plot the dimensionless pressure $p = P/\mu$, where $P$ being the real pressure.
}

\begin{figure}[h]
\centering
\begin{subfigure}[b]{0.33\linewidth}
	\centering
 	 \includegraphics[width=\linewidth]{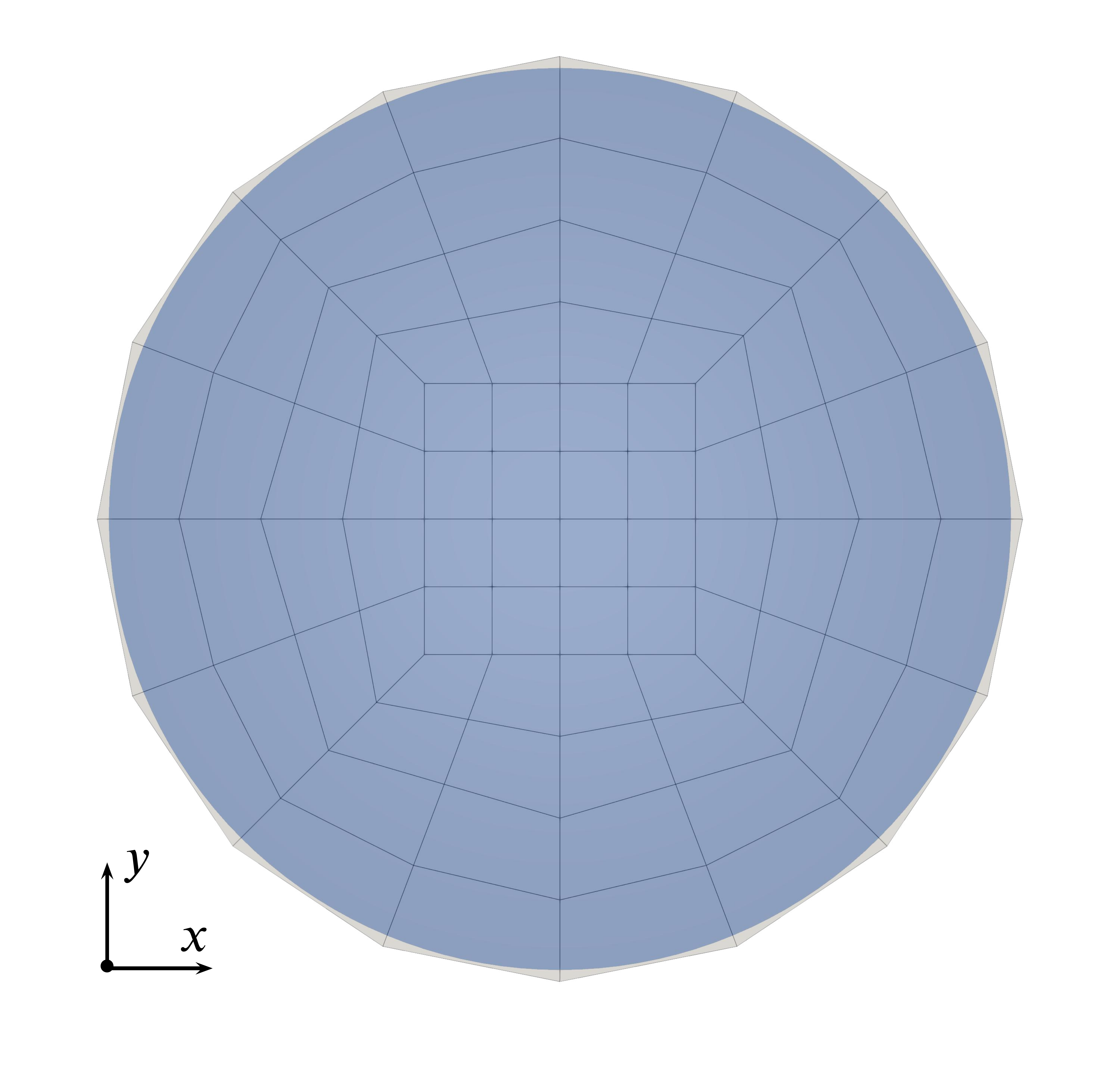}
	\caption{}
	\label{fig:plate_mesh}
\end{subfigure}
\begin{subfigure}[b]{0.5\linewidth}
	\centering
 	 \includegraphics[width=\linewidth]{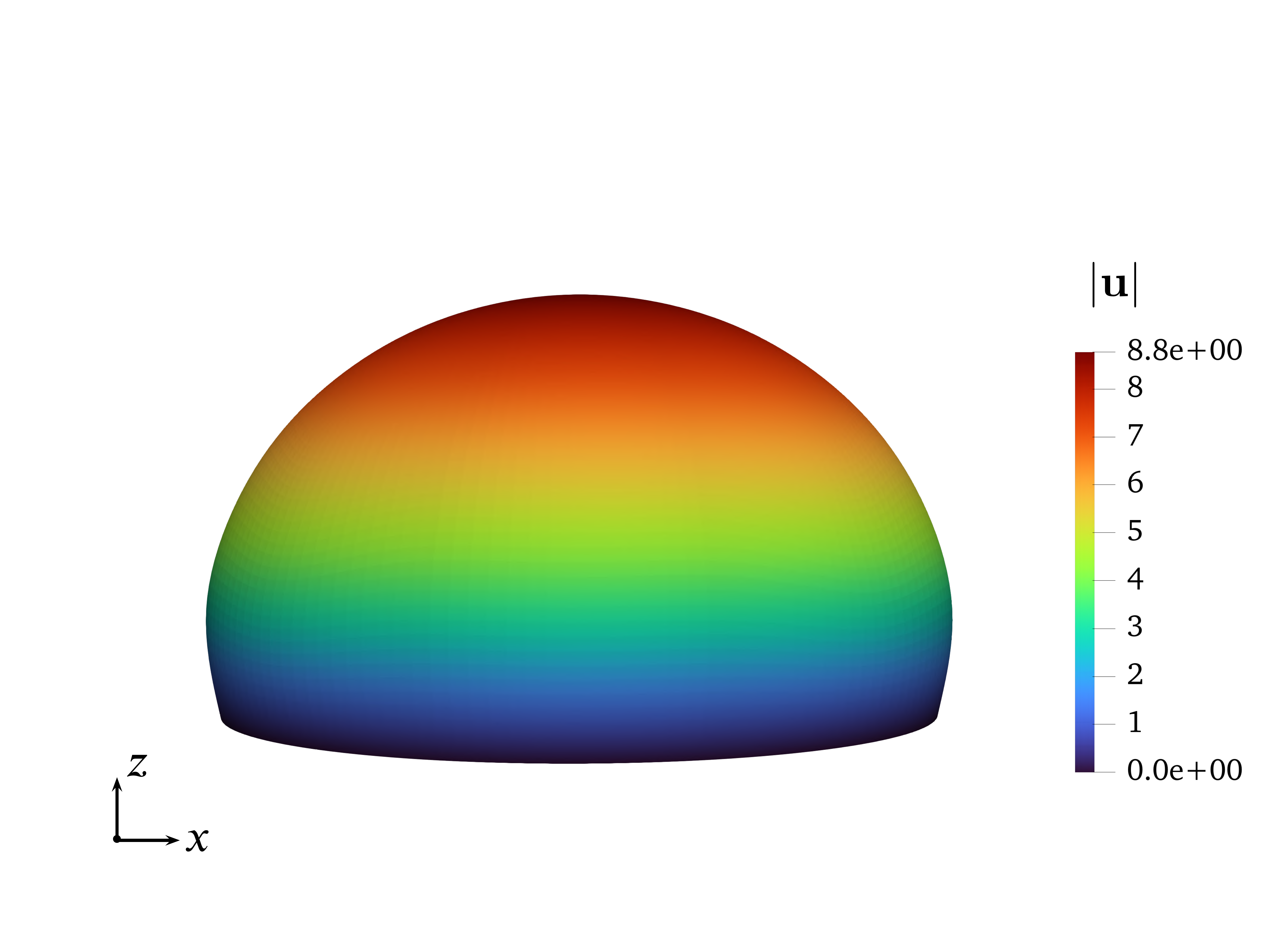}
	\caption{}
	\label{fig:plate_disp}
\end{subfigure}
\caption{(a) The initial control grid and the limit surface of the circular plate. (b) The deformed shape of the circular plate when dimensionless inflating pressure $p = 35$.}
\label{fig:circular_plate}
\end{figure}

The converged numerical result is obtained after two uniform refinements ($1440$ elements) and it was compared with the literature in \autoref{fig:circle_data}. The result has a slight offset compared to~\citet{hughes1983nonlinear} and ~\citet{Cirak:2001aa}, and perfectly agree with the recent work presented by~\citet{nama2020nonlinear}, whose shell theory is similar to~\cite{kiendl2015isogeometric} and was applied in biomechanics. This numerical example demonstrates the ability of the proposed method to solve a problem combining both membrane and bending deformations.
\begin{figure}
\centering
  \includegraphics[width=0.6\linewidth]{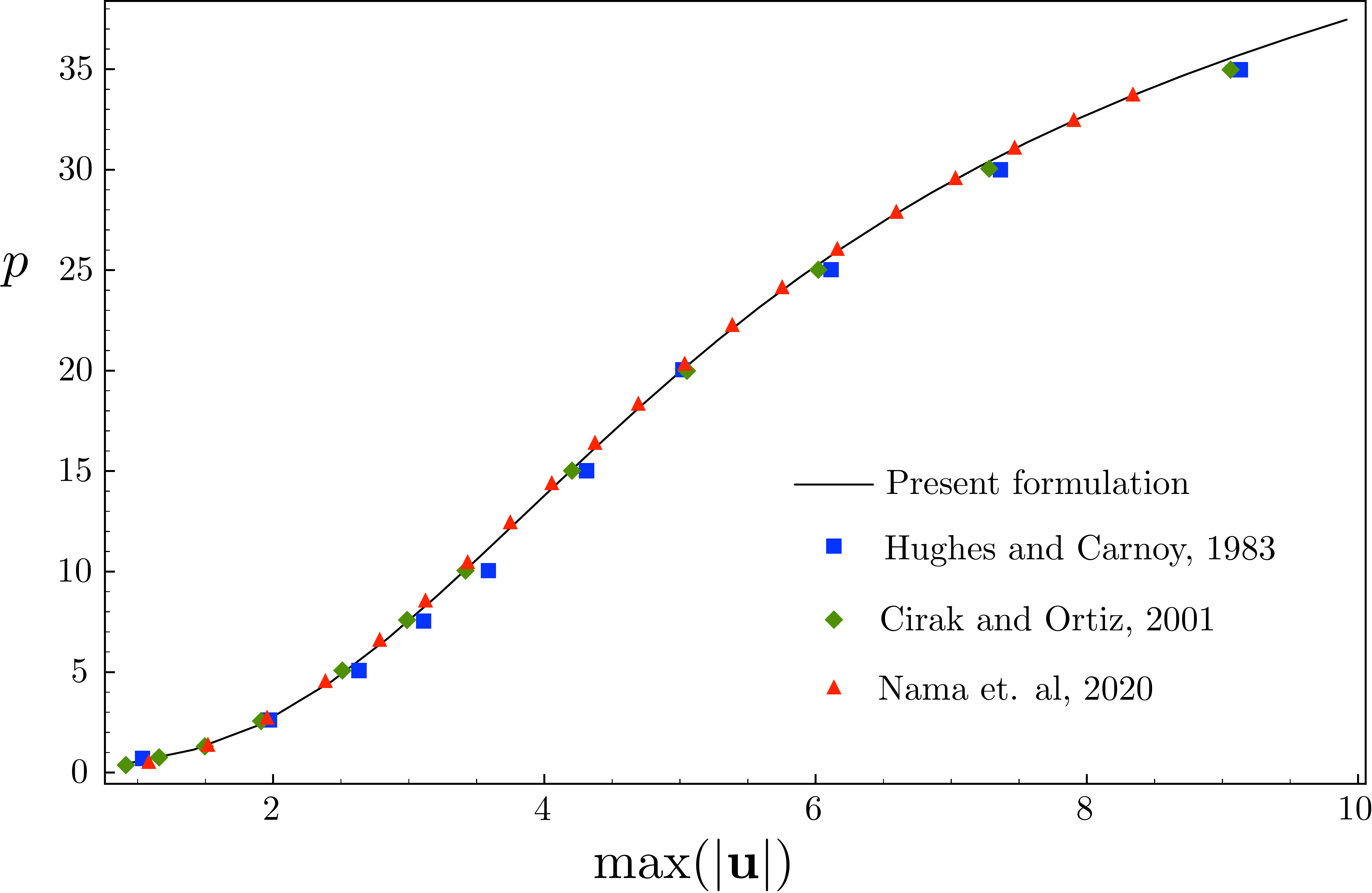}
\caption{Variation of the inflating pressure $p$ with the maximum displacement (max($|\mathbf{u}|$)) for inflation of a circular plate. Comparison
of the solution based on the present formulation with that presented in the literature.}
\label{fig:circle_data}
\end{figure}
\subsection{Inflation of a spherical balloon}
\label{sec:balloon}
The second example considered is the inflation of a spherical balloon, which was analytically studied in Section 6.3 of~\cite{holzapfel2002nonlinear}. 
It is also a widely used benchmark problem for hyperelatic shell formulations~\cite{Cirak:2001aa,chen2014explicit,kiendl2015isogeometric}. \autoref{fig:sphere_configs} shows the geometry for the spherical balloon. The balloon is inflated with a uniform pressure applied from the inside. For an incompressible Mooney-Rivlin model, the analytical solution of the internal pressure of the balloon is
\begin{equation}
    p = \frac{ 4 \bar h}{\bar R} \left[c_1 [\lambda^{-1} - \lambda^{-7}] - c_2[\lambda^{-5} - \lambda]\right],
\end{equation}
where $\bar R$ is the radius of the spherical shell in the reference configuration and $\bar h$ is the undeformed thickness.
We take the values $\bar R = 10 $ and $\bar h = 0.1$ for our simulations.
For the spherical hyperelastic shell, the stretching of the mid-surface is the same in all directions, thus $ \lambda_1 = \lambda_2 = \lambda$. 
Due to the incompressiblility constraint, the thickness stretch can be computed as $\lambda_3 = 1/\sqrt{\lambda}$. 
\begin{figure}
\centering
\begin{subfigure}[b]{0.5\linewidth}
	\centering
 	 \includegraphics[width=\linewidth]{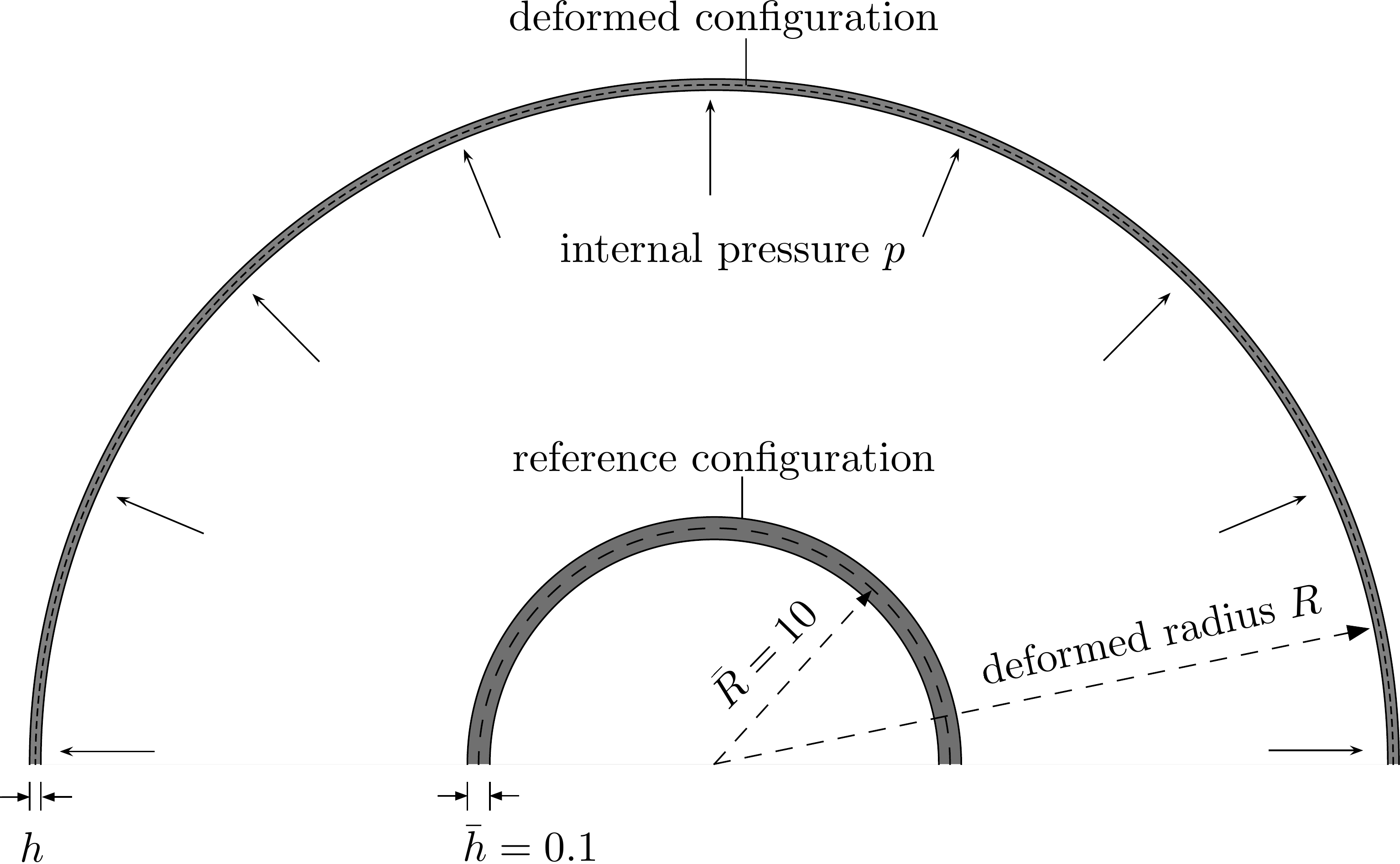}
	\caption{}
	\label{fig:sphere_configs}
\end{subfigure}
\begin{subfigure}[b]{0.45\linewidth}
	\centering
 	 \includegraphics[width=0.9\linewidth]{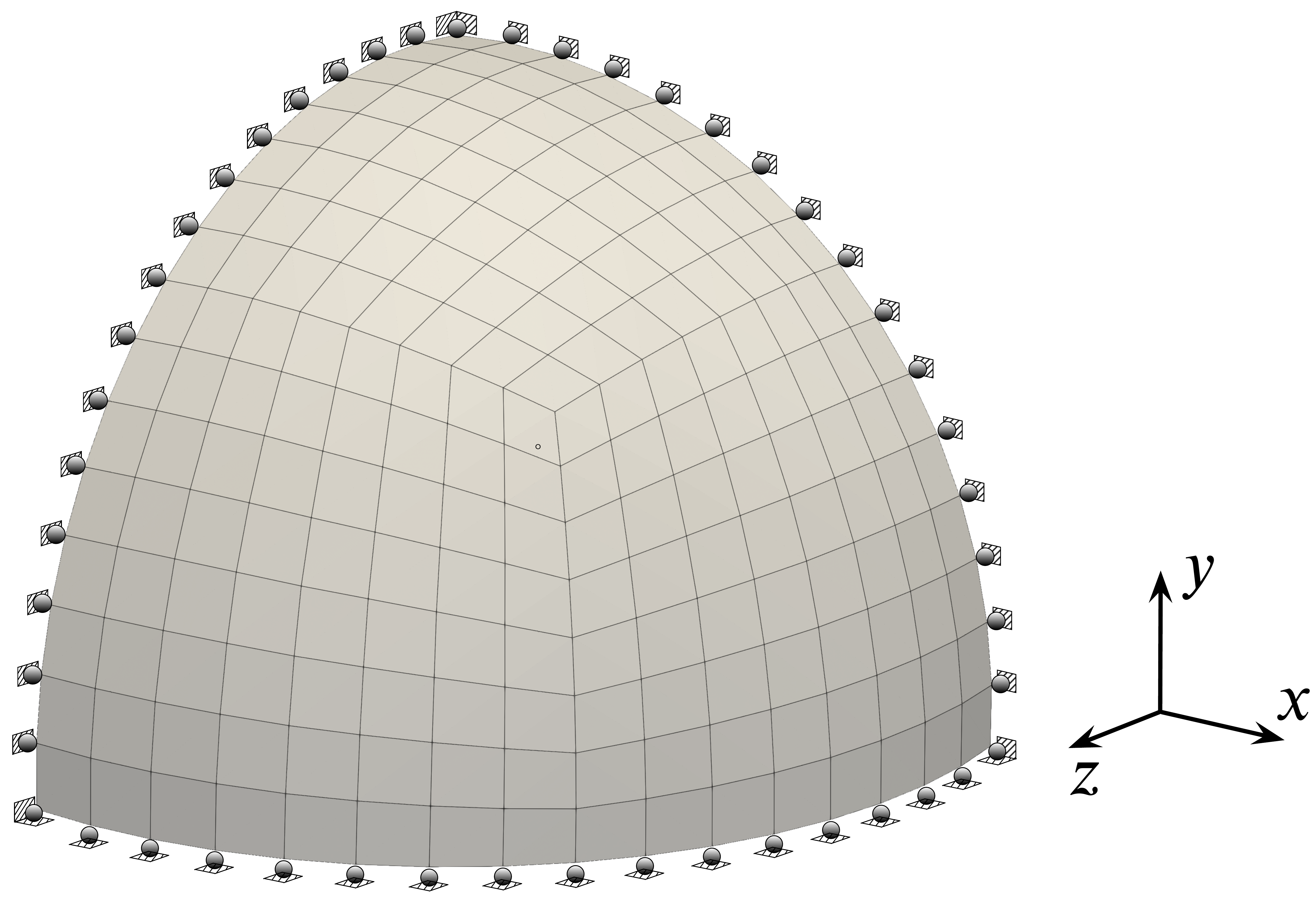}
	\caption{}
	\label{fig:sphere_mesh}
\end{subfigure}
\caption{(a) Geometrical setting of a spherical balloon in the reference and deformed configurations, showing only half of the cross section. (b) The mesh of a quarter of the hemisphere (one eighth of the full sphere) is used for the numerical test and the symmetry assumption is adopted by constraining the corresponding degree of freedom for the edges shown in the figure.}
\label{fig:spheres}
\end{figure}

\autoref{fig:sphere_mesh} shows the control mesh with $192$ elements of one quarter of the hemisphere used for the numerical simulation. Symmetry boundary conditions are applied to three edges by constraining the corresponding degrees of freedom as shown. Two sets of material parameters are tested for this example. The parameters for the first case are $c_1 = 0.5 \mu$ and $c_2 = 0$, where $\mu = 4.225 \times 10^{5}$. 
Because $c_2$ is set to zero, the material model is reduced to an incompressible neo-Hookean. The second case of parameters is $c_1 = 0.4375\mu$ and $c_2 = 0.0625\mu$, thus $c_1/c_2 = 7$.  \autoref{fig:sphere_curves} shows that the numerical results for both tests perfectly agree with the analytical solutions.
\begin{figure}
\centering
  \includegraphics[width=\linewidth]{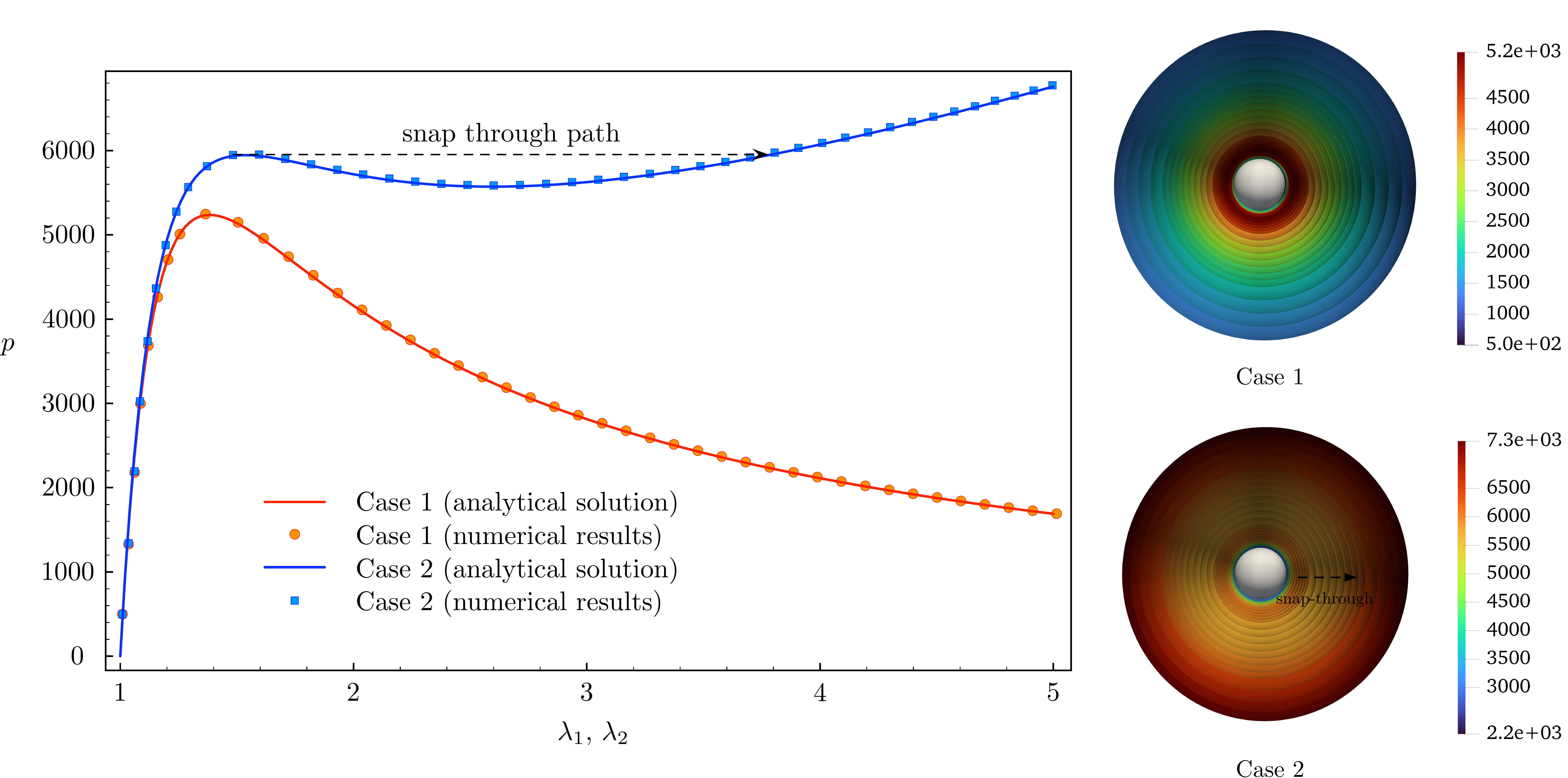}
\caption{{Variation of the inflation pressure against the in-plane stretch values ($\lambda_1, \lambda_2$).} The numerical and analytical stretch - pressure curves for the inflated balloons (Case 1: neo-Hookean model,  $c_1 = 0.5 \mu$ and $c_2 = 0$; Case 2: Mooney--Rivlin model, $c_1 = 0.4375\mu$ and $c_2 = 0.0625\mu$ , where $\mu = 4.225 \times 10^{5}$ ). Excellent agreement is observed between the analytical and numerical solutions for both the material models. The inflated profiles of two cases are shown to the right. A full undeformed balloon is shown uncoloured for both cases and the deformed profiles for increasing load steps are shown as hemispheres. The colour indicates the value of internal pressure. For the second case, the balloon may experience snap-through during inflation. }
\label{fig:sphere_curves}
\end{figure}
For the first case, the internal pressure of the spherical balloon reaches the limit point when $\lambda = 1.38$. If the volume of the balloon keeps increasing, the internal pressure gradually decreases. However, in the second case, the internal pressure will first decrease after the limit point and then rise again, exceeding the pressure at the limit point. The complex nonlinear response of the structures are captured by the arc-length method. \autoref{fig:sphere_curves} also shows the inflated profiles of the spherical balloon for both cases and the internal pressures values are indicated with the colouring. 
If the two balloons are both inflated using pressure-control, the neo-Hookean balloon will explode when it is pressurised beyond the limit point, because that is the maximum pressure it can withstand.
The Mooney--Rivlin balloon may undergo a sudden and large deformation at the limit point due to the increased pressure and reach a new equilibrium state. 
This is the snap-through phenomenon which is normally considered as an instability of the structure.
With the nonlinear algorithm and arc-length method, the proposed method can easily capture this phenomenon numerically even for complex structures.


\subsection{Loss of symmetry during the inflation of a toroidal thin shell}
This example considers the bifurcation of an inflated toroidal thin shell. A torus is the simplest example of a genus 1 orientable surface. Toroidal membranes and shells are widely applied in engineering applications such as tyres, air springs, soft grippers and inflated actuators. The bifurcation instability (loss of symmetry) of a toroidal membrane is examined semi-analytically in~\cite{pamulaparthi2019instabilities}. Here, a similar problem is analysed numerically for a thin shell. 
\begin{figure}
\centering
\begin{subfigure}[b]{0.4\linewidth}
	\centering
 	 \includegraphics[width=\linewidth]{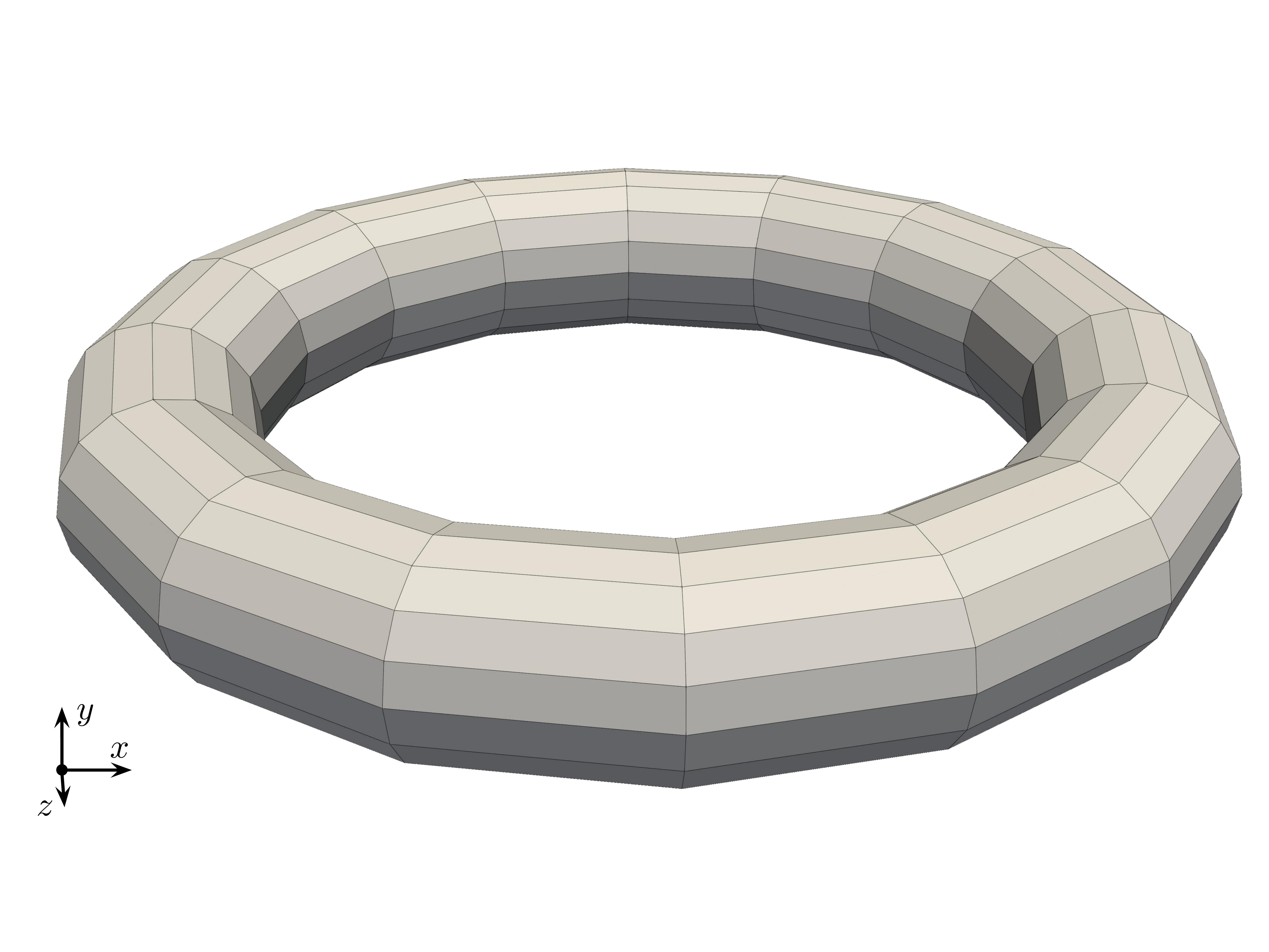}
	\caption{}
	\label{fig:torus_mesh}
\end{subfigure}
\begin{subfigure}[b]{0.4\linewidth}
	\centering
 	 \includegraphics[width=\linewidth]{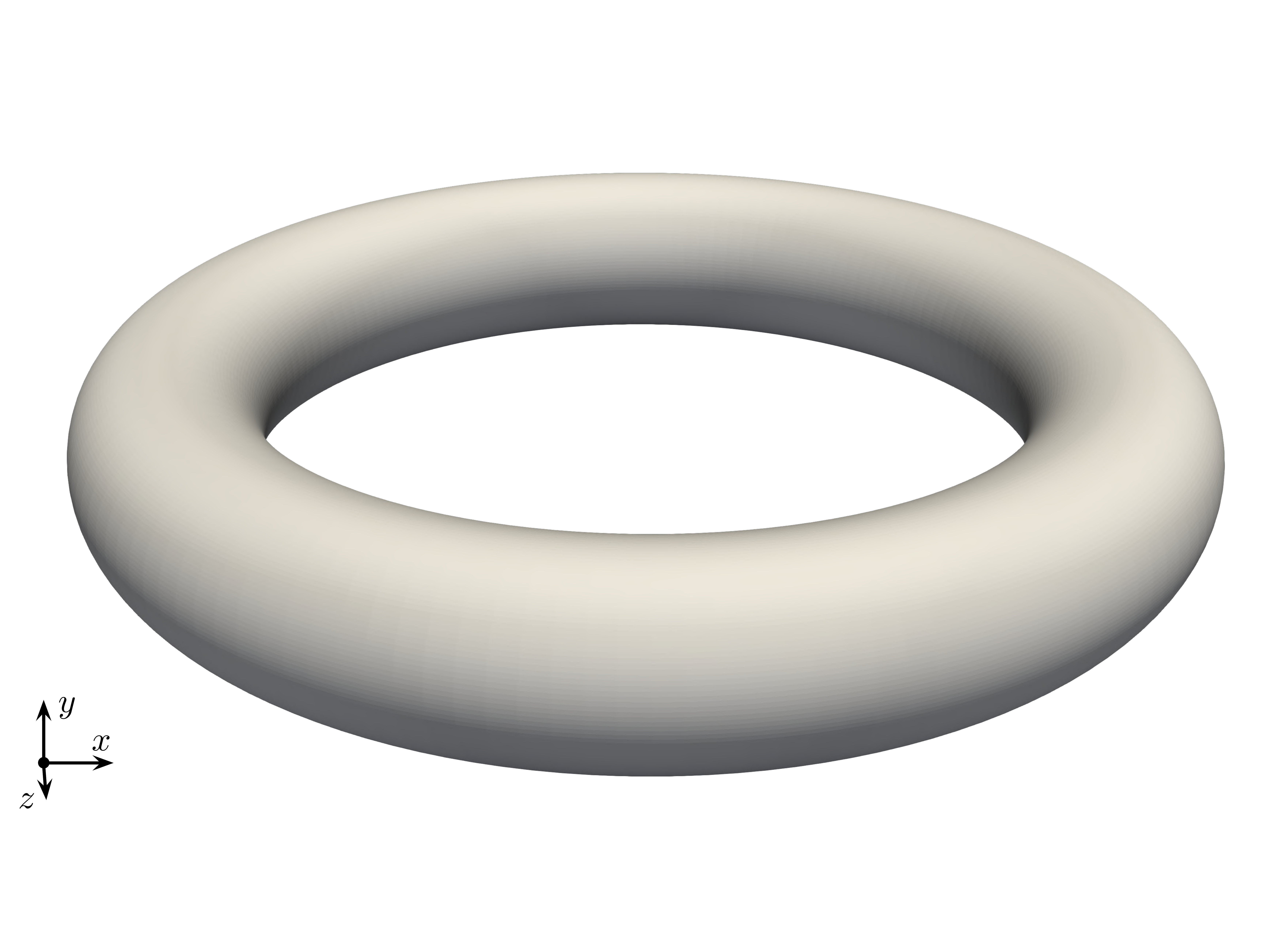}
	\caption{}
	\label{fig:torus_limit}
\end{subfigure}
\caption{A smooth torodial surface is modelled with a relatively coarse control grid. (a) The control grid with $256$ elements. (b) The limit surface of the toroidal thin shell.}
\label{fig:torus_inflation_1}
\end{figure}
The structure is modelled as an enclosed hyperelastic thin shell and the inflation is simulated. The internal pressure of the thin shell is applied incrementally as an external force and the displacements are solved by using the Newton-Raphson method for each load step. After the deformed equilibrium state is achieved for each load step, an eigenvalue analysis of the stiffness matrix, as introduced in Section~\ref{sec:eigenvalue}, is performed to check stability. When the stiffness matrix has zero eigenvalues, the structure is in an unstable state where sudden geometric changes may occur. The eigenvector corresponding to the zero eigenvalue indicates the direction of the possible sudden change. If it breaks the original symmetry of the structure, the inflated toroidal thin shell will bifurcate to a new non-symmetric branch. A control grid with $256$ elements shown in \autoref{fig:torus_mesh} renders a limit mid-surface of the toroidal thin shell shown in \autoref{fig:torus_limit}. The minimum bounding box for the limit surface are $[-10.6522,10.6522]\times[-1.80474,1.80474]\times[-10.6522,10.6522]$. It is a relatively slender torus and its aspect ratio is approximately $4.9$. The thickness of the shell is set to $0.01 $ and the same material parameters $c_1$ and $c_2$ as for the spherical balloon example in Section~\ref{sec:balloon}  are adopted here. 

\autoref{fig:torus_plot} shows the evolution of the internal pressure against the enclosed volume for the toroidal thin shell modelled with Catmull-Clark subdivision surfaces.
When the thin shell is inflated to the limit point, the stiffness matrix has a pair of zero eigenvalues whose corresponding eigenmodes are orthogonal and they are both non-symmetric.
One of the eigenmodes is shown in the figure. The limit point is the earliest point at which bifurcation is likely to occur. The eigenmode is used to perturb the structure in the following load steps and a bifurcated solution is branched off from the principal solution. The inflated shapes of the toroidal thin shell for both branches are selectively plotted in \autoref{fig:torus_inflation_2}. 
\begin{figure}
\centering
\begin{subfigure}[b]{0.6\linewidth}
	\centering
 	 \includegraphics[width=\linewidth]{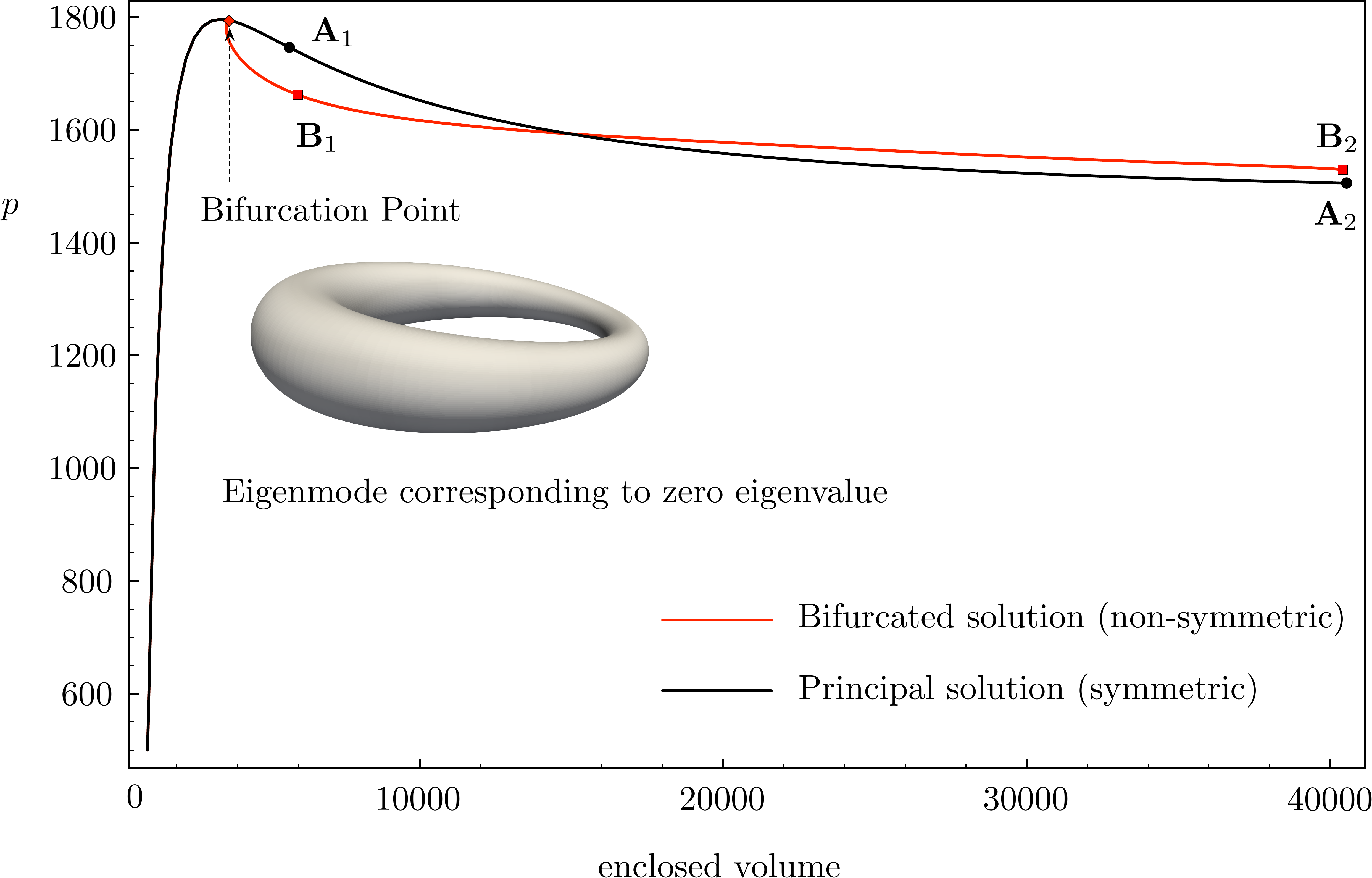}
	\caption{The variation of the applied pressure against the enclosed volume.}
	\label{fig:torus_plot}
\end{subfigure}
\begin{subfigure}[b]{0.45\linewidth}
	\centering
 	 \includegraphics[width=\linewidth]{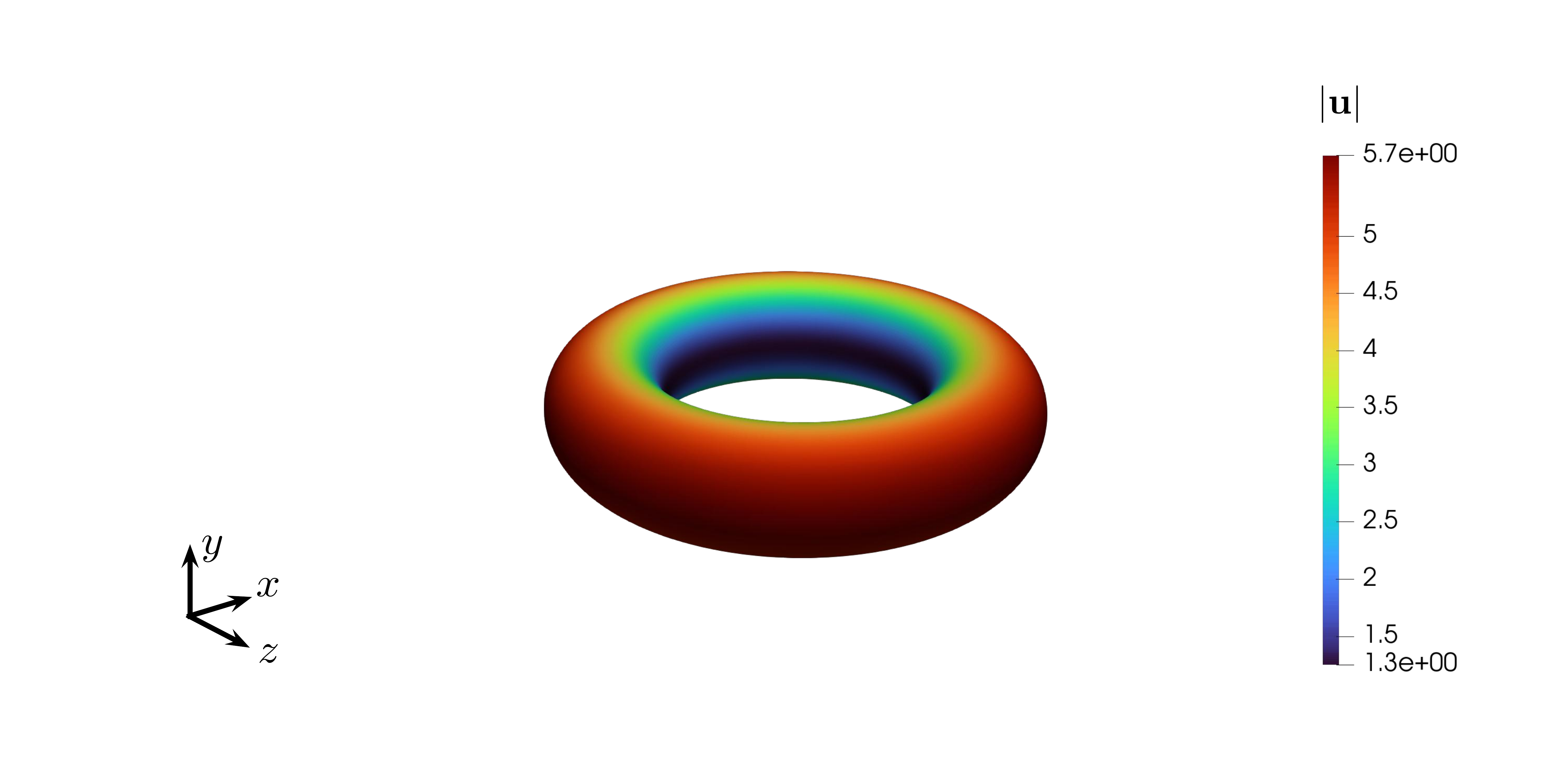}
	\caption{$|\mathbf u|$ at $\mathbf A_1$}
	\label{fig:torus_p15}
\end{subfigure}
\begin{subfigure}[b]{0.45\linewidth}
	\centering
 	 \includegraphics[width=\linewidth]{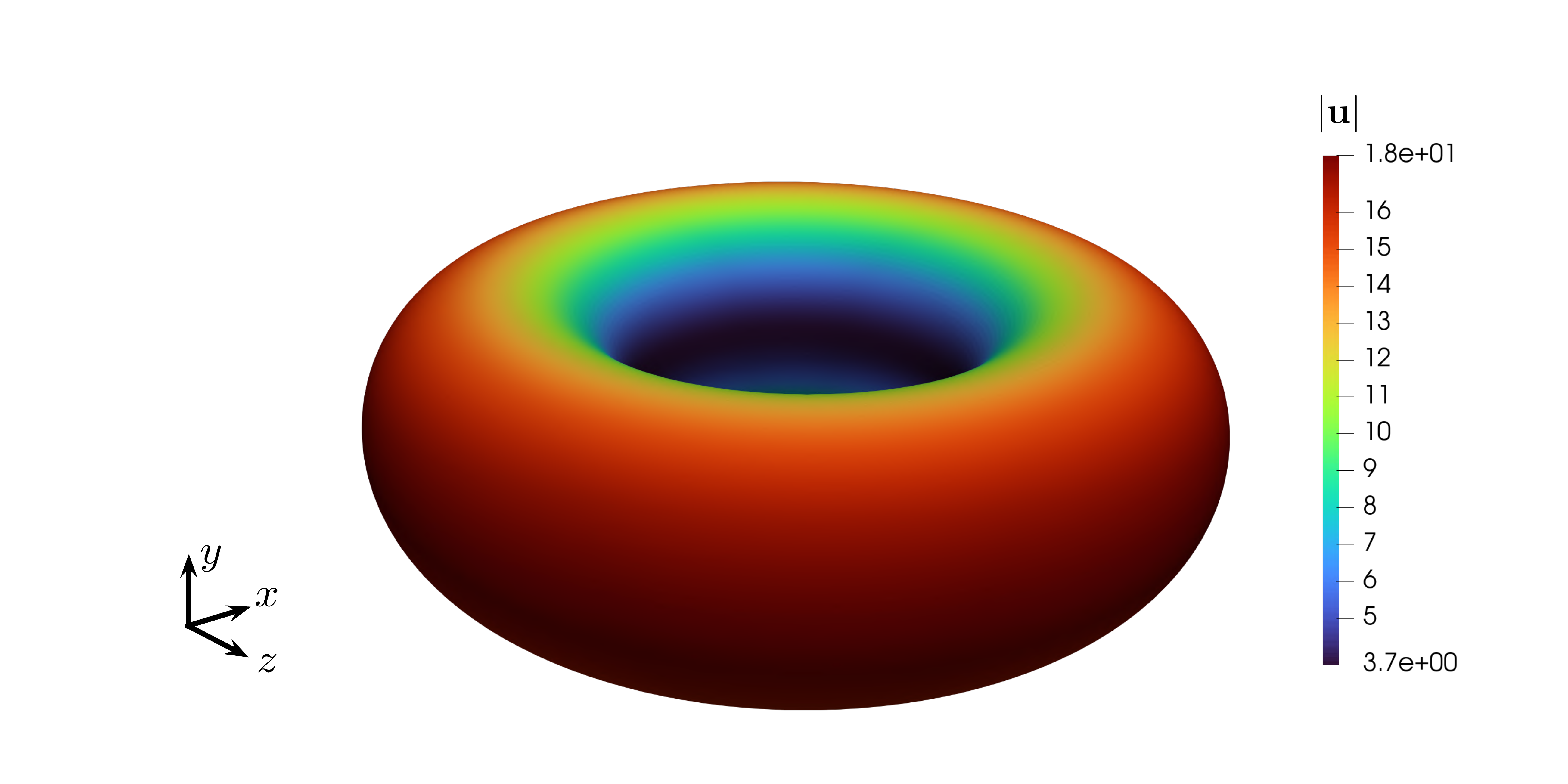}
	\caption{$|\mathbf u|$ at $\mathbf A_2$}
	\label{fig:torus_p47}
\end{subfigure}
\begin{subfigure}[b]{0.45\linewidth}
	\centering
 	 \includegraphics[width=\linewidth]{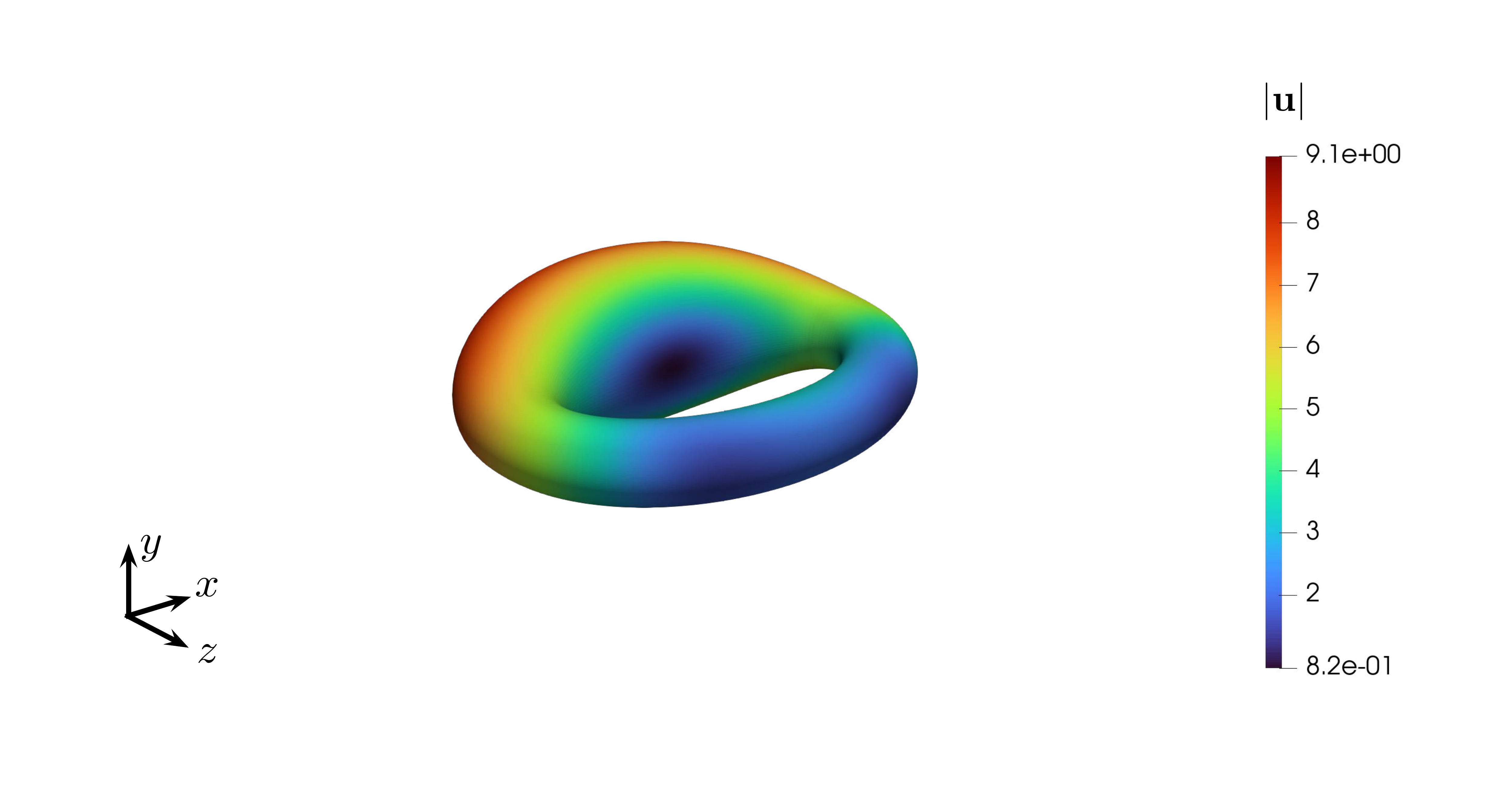}
	\caption{$|\mathbf u|$ at $\mathbf B_1$}
	\label{fig:torus_b24}
\end{subfigure}
\begin{subfigure}[b]{0.45\linewidth}
	\centering
 	 \includegraphics[width=\linewidth]{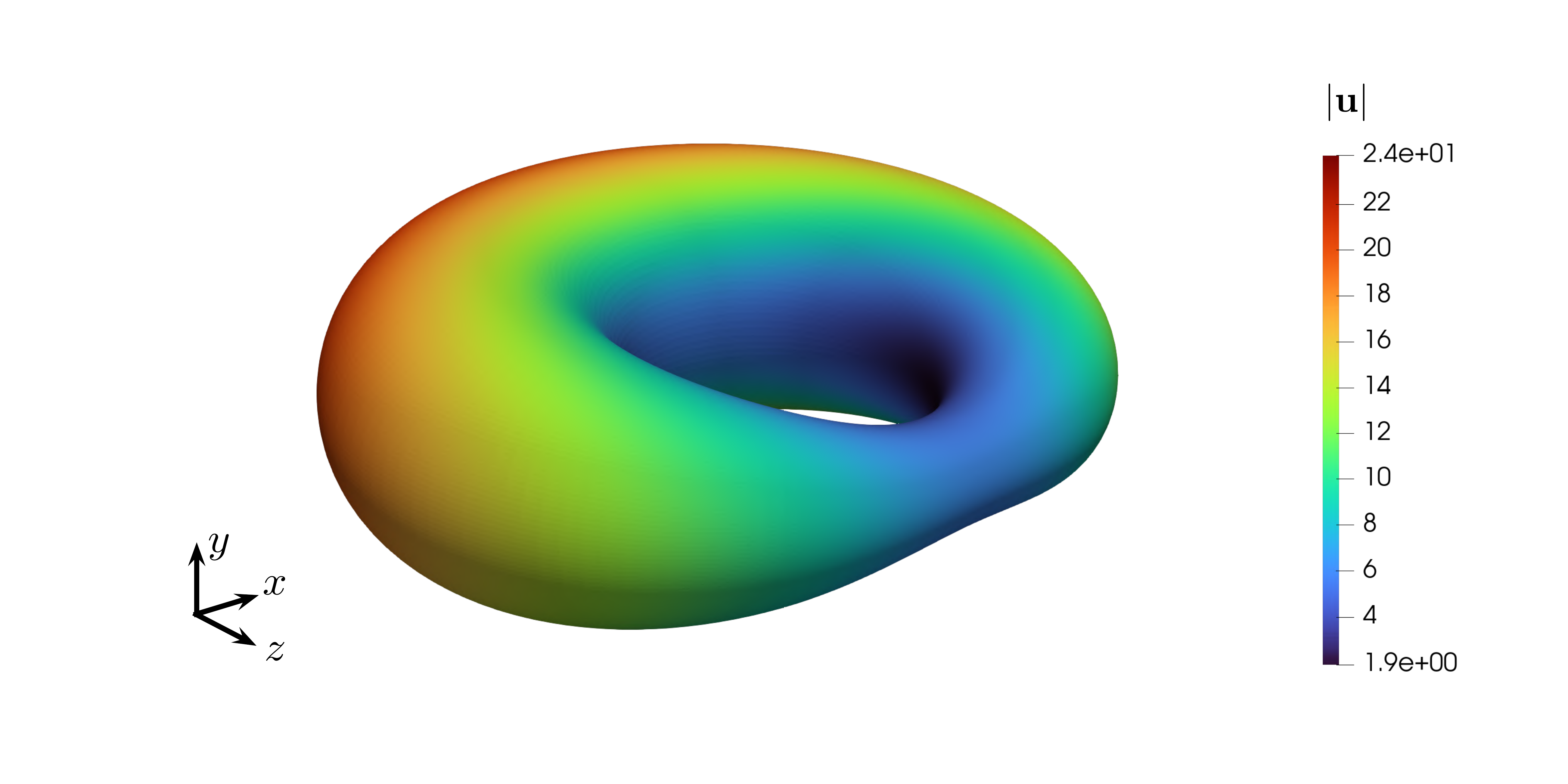}
	\caption{$|\mathbf u|$ at $\mathbf B_2$}
	\label{fig:torus_b73}
\end{subfigure}
\caption{Numerical computation of the variation of the applied pressure against the enclosed volume for the inflated toroidal thin shell is plotted in (a). 
For this case, bifurcation occurs immediately after the solution reaches the limit point.
The asymmetric eigenmode at the bifurcation point is also shown. Deformed shapes of the toroidal thin shell during inflation for both principal ((b) and (c)) and bifurcated ((d) and (e)) solutions are shown. The locations of the deformation states $\mathbf A_1$, $\mathbf A_2$, $\mathbf B_1$, $\mathbf B_2$ in the pressure-volume curves are indicated in \autoref{fig:torus_plot}.}
\label{fig:torus_inflation_2}
\end{figure}
\subsection{Inflation of an elastic airbag with folds and wrinkles}
The final example considers the problem of inflating an elastic airbag. This a classic example to study wrinkles of membranes and thin shells in the literature~\cite{lu2001finite,Cirak:2001aa,chen2014explicit,nakashino2020geometrically}, where the inflation of a airbag is considered as a highly nonlinear problem and the solution of the problem is non-unique and challenging to compute.

A subdivision surface with $256\,(16\times 16)$ elements is generated using a control mesh for a square plate in order to model a half airbag. Symmetry boundary conditions are applied to the four edges by constraining their degrees of freedom corresponding to the displacements in $z\mbox{-}$direction and the in-plane rigid body motion is also eliminated. The constitutive relation for this problem is degenerated to the Saint Venant-Kirchhoff model since the airbag textile is often considered as inextensible and is easy to wrinkle and fold. The energy density is expressed as
\begin{equation}
  W(\mathbf{E})=\frac{\mathit{E} \, \nu}{ 2[1+\nu][1-2\nu]}[\operatorname{tr}( \mathbf{E})]^2+\frac{\mathit{E} }{ 2[1+\nu]}\mathbf{E}: \mathbf{E}.
\end{equation}
The Young's modulus and Poisson's ratio are set as $\mathit{E} = 5\times 10^8  \text{ and }\nu = 0.4$, respectively. The length of the side of the square airbag is $1$ and its thickness is set to $0.001$. The airbag is inflated by incrementally increasing the pressure until the inflating pressure reaches the value $p = 5000$. 
Due to the highly nonlinear nature of the problem, we observe a bifurcation in the solution very close to the reference configuration.
As a result, the initial load step determines the solution branch.
\autoref{fig:airbag_200_u} shows the final deformed shape of the airbag for the case 
when the internal pressure $\Delta p = 200$ is selected as the first load step. 
The complex folds and wrinkles seen in \autoref{fig:airbag_200_u} are easily captured with our thin shell formulation. 

\begin{figure}
\centering
\begin{subfigure}[b]{0.6\linewidth}
	\centering
 	 \includegraphics[width=\linewidth]{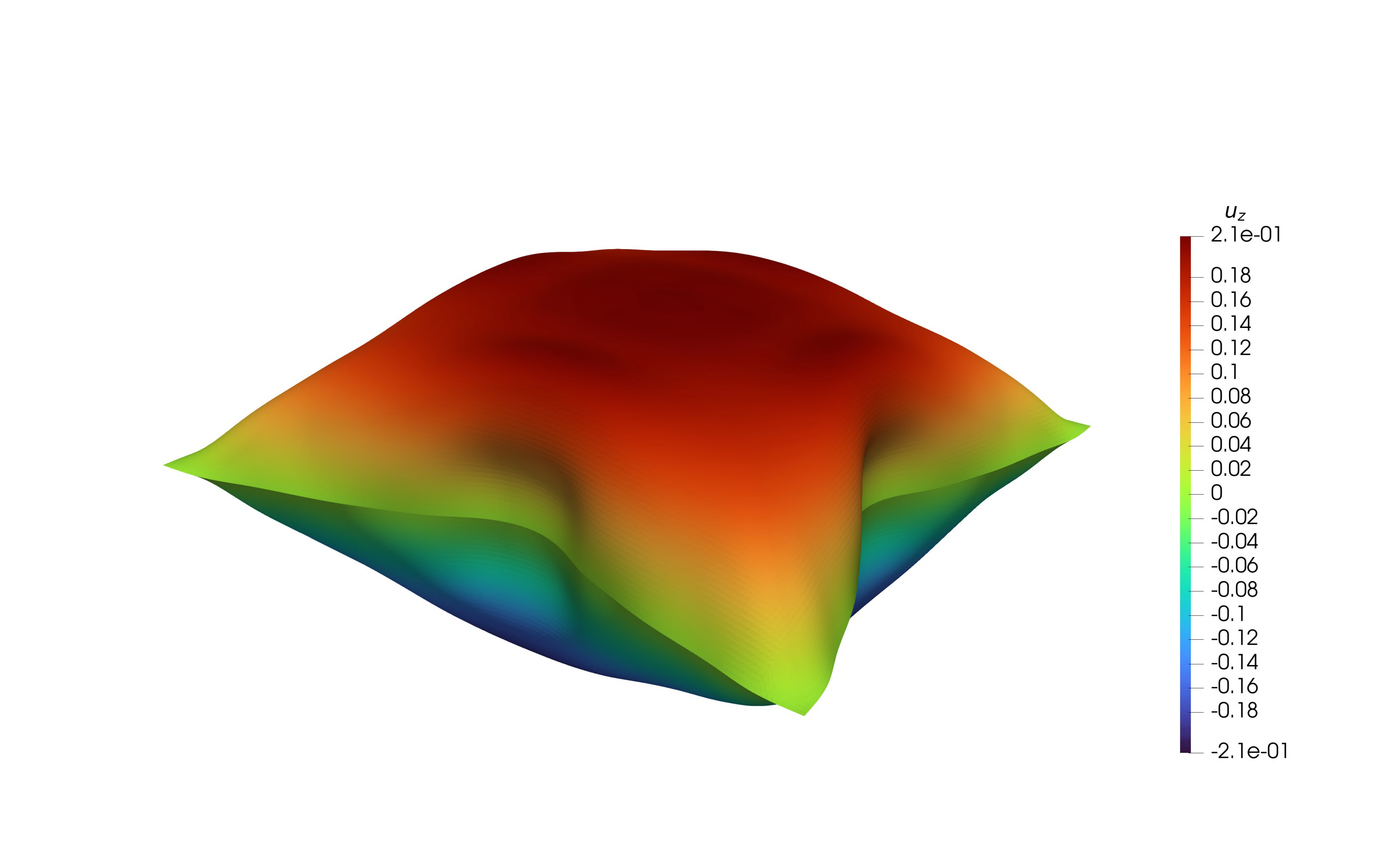}
	\caption{Vertical displacement $u_z$ of the mid-surface at $p=5000$.}
	\label{fig:airbag_200_u}
\end{subfigure}
\begin{subfigure}[b]{0.6\linewidth}
	\centering
 	 \includegraphics[width=\linewidth]{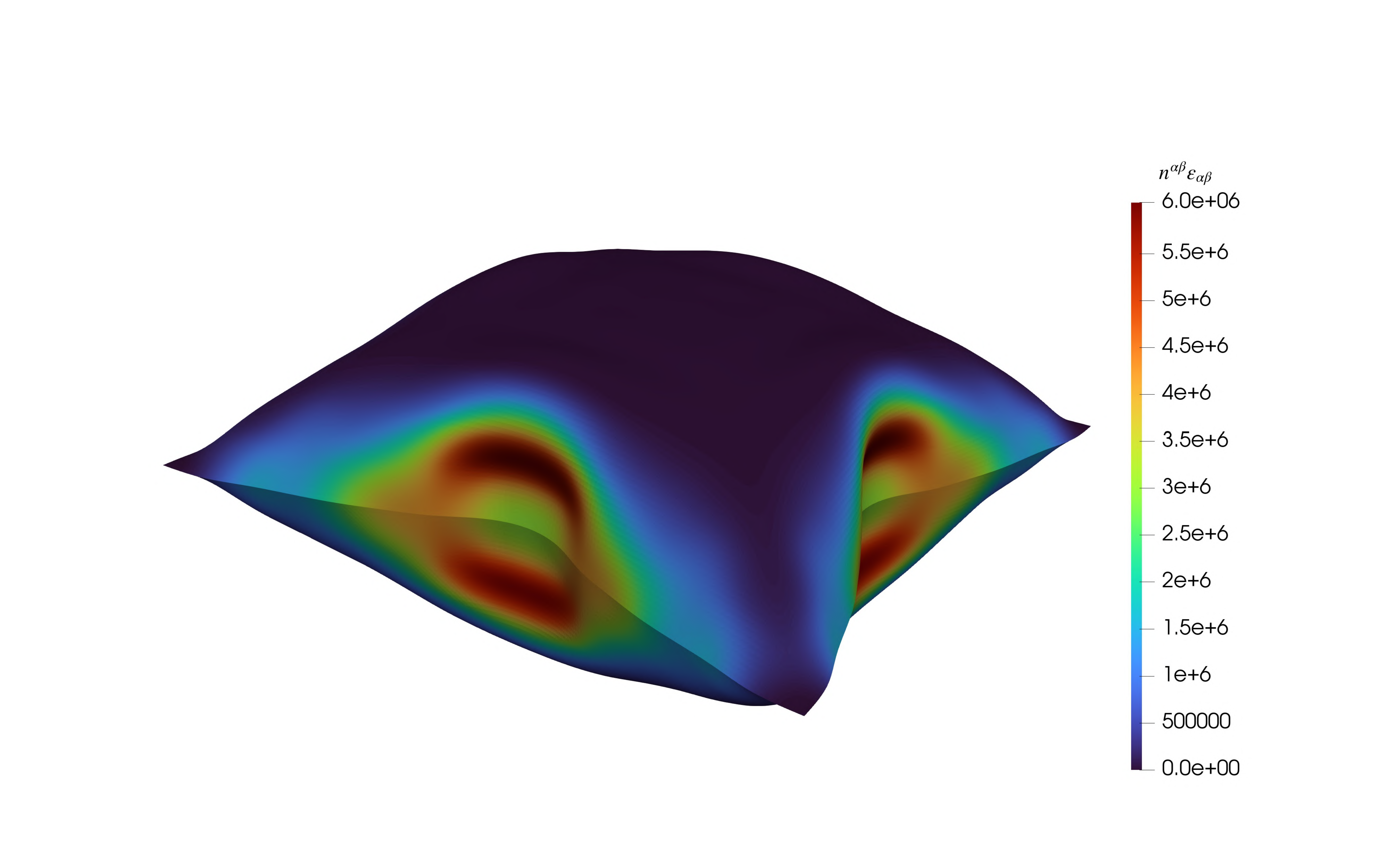}
	\caption{Energy density $n^{\alpha\beta}\varepsilon_{\alpha\beta}$ at the mid-surface at $p=5000$.}
	\label{fig:airbag_200_E}
\end{subfigure}
\caption{Results for the simulation of inflation of a square airbag modelled using a St. Venant--Kirchhoff constitutive model at an internal pressure of $p=5000$. 
}
\label{fig:airbags}
\end{figure}
\section{Summary and Conclusions}
An isogeometric approach for the analysis of inflated hyperelastic thin shell structures has been proposed. The Kirchhoff-Love hypothesis is followed in the thin shell formulation. 
The incompressible Mooney-Rivlin model is adopted to describe the material of the inflated hyperelastic thin shells. 
Based on the principle of virtual work, the weak form of governing equations of thin shell structures is formulated. 
The complex nonlinear response of the inflated hyperelastic thin shell are numerically simulated by the arc-length method, where the displacements and internal pressure are solved incrementally. For each increment, the displacements and internal pressure is solved using the Newton-Raphson method. Both the geometry and deformation field are discretised using the Catmull-Clark subdivision bases and a finite element framework with $C^1$ continuity is established. Two types of the global instabilities of the inflated thin shell structures can be simulated, which are the snap-through and bifurcation buckling. The proposed method is first validated by taking an inflated circular plate as an example, which has been widely studied in the literature. A numerical example of an inflated spherical shell further verifies the formulation analytically. Moreover, it shows the ability of the proposed method to reveal the snap-through phenomenon of hyperelastic shells with enclosed volumes.
Thereafter, a torodial hyperelastic thin shell is simulated  to investigate the bifurcation of solution from a principal symmetric mode.
After solving the displacements for each inflation increment, the eigenvalues of the stiffness matrix have been checked, and the structure is perturbed with the eigenvector if the corresponding eigenvalue is zero, thus inducing bifurcation. The numerical simulation shows that the torodial thin shell analysed in the present work may lose its symmetry when the inflation reaches its limit point. Finally, the inflation  of an airbag is simulated, which has demonstrated the ability of proposed method to capture complex states such as wrinkles and folds in finitely deformed thin shells. 

For inflation problems, the membrane stiffness of a highly deformable hyperelastic thin shell dominates the response.
Due to strong kinematic and constitutive nonlinearities, the mechanical response of a thin shell can be challenging to simulate - especially close to bifurcation points.
Care must be taken when selecting the initial loading values, the arc-length parameters, and the scaling factor used in the eigenvalue analysis.
The focus here is however on the formulation and the demonstration of its capabilities.

\section*{Acknowledgements}
This work was supported by the UK Engineering and Physical Sciences Research Council (EPSRC) grants EP/R008531/1 and EP/V030833/1.
We also thank for the support from the Royal Society International Exchange Scheme IES/R1/201122.
Paul Steinmann gratefully acknowledges financial support for this work by the Deutsche Forschungsgemeinschaft under GRK2495, projects B \& C. 
Tiantang Yu acknowledges the support from the National Natural Science Foundation of China (NSFC) under Grant Nos.11972146 and 12272124.
\appendix
\section{Appendix}
\subsection{Variations of normal vector and thickness stretch}
\label{sec:Appendix_1}
To simplify the expression, one denotes the normal vector as
\begin{equation}
\mathbf{a}_3 = J^{-1} \tilde{\mathbf{a}}_3.
\end{equation}
where
\begin{equation}
\tilde{\mathbf{a}}_3 = \mathbf{a}_{1}\times\mathbf{a}_{2}.
\end{equation}
Then, the first and second variations of the normal vector are computed as
\begin{equation}
\begin{aligned}
\delta_{r}\mathbf{a}_{3}&=J^{-1}\delta_{r}\tilde{\mathbf{a}}_{3}-J^{-2}[\delta_{r}J] \tilde{\mathbf{a}}_{3},\\
\delta_{s}\delta_{r}\mathbf{a}_{3} &= J^{-1} [\delta_{s}\delta_{r}\tilde{\mathbf{a}}_{3}]-J^{-2} [\delta_{s}\delta_{r}J] \tilde{\mathbf{a}}_{3}-J^{-2} [\delta_{r}J][\delta_{s} \tilde{\mathbf{a}}_{3}]-J^{-2} [\delta_{s}J][ \delta_{r}\tilde{\mathbf{a}}_{3}] + 2 J^{-3} [\delta_{r}J][ \delta_{s} J] \tilde{\mathbf{a}}_{3},
\end{aligned}
\end{equation}
where
\begin{equation}
\begin{aligned}
\delta_{r}\tilde{\mathbf{a}}_{3} &= \delta_{r}\mathbf a_1 \times \mathbf a_2 + \mathbf a_1 \times \delta_{r}\mathbf a_2,\\
\delta_{r}J&=\mathbf{a}_{3} \cdot \delta_{r}\tilde{\mathbf{a}}_{3},\\
\delta_{s}\delta_{r}\tilde{\mathbf{a}}_{3} & = \delta_{r}\mathbf a_1 \times \delta_{s}\mathbf a_2 + \delta_{s}\mathbf a_1 \times \delta_{r}\mathbf a_2, \\
\delta_{s}\delta_{r}J &= J^{-1}\left[\delta_{s}\delta_{r}\tilde{\mathbf{a}}_{3} \cdot \tilde{\mathbf{a}}_{3}+\delta_{r}\tilde{\mathbf{a}}_{3} \cdot \delta_{s}\tilde{\mathbf{a}}_{3}-[\delta_{r}\tilde{\mathbf{a}}_{3} \cdot \mathbf{a}_{3}][\delta_{s}\tilde{\mathbf{a}}_{3} \cdot \mathbf{a}_{3}]\right].
\label{eq:variations_a3}
\end{aligned}
\end{equation}
Due to the incompressibility constraint of the material, the thickness stretch can be expressed as
\begin{equation}
\lambda_3 = \frac{|\bar{\mathbf a}_1 \times \bar{\mathbf a}_2|}{|{\mathbf a}_1 \times {\mathbf a}_2|} = \bar{J} J^{-1} ,
\end{equation}
whose first and second variations are expressed as
\begin{equation}
\begin{aligned}
\delta_r \lambda_3 
=& -\bar{J} [ J]^{-2}\delta_r J,\\
\delta_s \delta_r \lambda_3 
 =& -\bar{J} \left[ -2[J]^{-3}\delta_sJ \, \delta_r J + [J]^{-2} \delta_s\delta_r J\right].\\
\end{aligned}
\end{equation}
\end{nolinenumbers}
\bibliographystyle{cas-model2-names}      
\bibliography{cas-refs}
\end{document}